\documentclass[]{aspm}
\usepackage{epsfig}
\usepackage{latexsym}
\usepackage[arrow,matrix,graph,frame,poly,arc,tips]{xy}
\usepackage{amscd,amssymb,latexsym,subfigure,verbatim,hyperref,epsfig,amsmath,su
pertabular}
\xyoption{all}
\usepackage{calc}
\articleinfo{}{}{}

\renewcommand{\P}{{\mathbb{P}}}
\newcommand{\K}{{\mathbb{K}}}
\newcommand{\Z}{{\mathbb{Z}}}
\newcommand{\C}{{\mathbb{C}}}
\newcommand{\N}{{\mathbb{N}}}
\newcommand{\LL}{{\mathbb{L}}}
\newcommand{\Q}{{\mathbb{Q}}}

\newcommand{\A}{{\mathcal{A}}}
\newcommand{\CL}{{\mathcal{C}}}
\newcommand{\Dp}{{D^p(\A)}}

\newtheorem{defn0}{Definition}
\newtheorem{prop0}[defn0]{Proposition}
\newtheorem{conj0}[defn0]{Conjecture}
\newtheorem{thm0}[defn0]{Theorem}
\newtheorem{lemma0}[defn0]{Lemma}
\newtheorem{corollary0}[defn0]{Corollary}
\newtheorem{example0}[defn0]{Example}
\newtheorem{problem0}[defn0]{Problem}

\newenvironment{defn}{\begin{defn0}}{\end{defn0}}

\newenvironment{conj}{\begin{conj0}}{\end{conj0}}
\newenvironment{thm}{\begin{thm0}}{\end{thm0}}

\newenvironment{corollary}{\begin{corollary0}}{\end{corollary0}}
\newenvironment{exm}{\begin{example0}\rm}{{$\Diamond$}\end{example0}}
\newenvironment{prob}{\begin{problem0}\rm}{\end{problem0}}

\newcommand{\Tor}{\mathop{\rm Tor}\nolimits}
\newcommand{\Ext}{\mathop{\rm Ext}\nolimits}
\newcommand{\im}{\mathop{\rm Im}\nolimits}

\newcommand{\codim}{\mathop{\rm codim}\nolimits}
\newcommand{\pdim}{\mathop{\rm pdim}\nolimits}

\newcommand{\XA}{{M_{\mathcal A}}}
\newcommand{\xA}{{X_{\mathcal A}}}
\newcommand{\g}{\mathfrak{g}}
\newcommand{\hA}{\mathfrak{h}_{\mathcal A}}

\newcommand{\ann}{\mathop{\rm ann}\nolimits}

\title[Hyperplane Arrangements: Computations and Conjectures]
{Hyperplane Arrangements: Computations and Conjectures}

\author[]{Hal Schenck}
\thanks{Schenck supported by  NSF 0707667}

\address{Schenck: Mathematics Department, University of Illinois Urbana-Champaign, Urbana IL 61801 USA}

\email{schenck@math.uiuc.edu}

\rcvdate{Month Day, Year}
\rvsdate{Month Day, Year}



\begin{document}

\begin{abstract}
This paper provides an overview of selected results and open problems in the 
theory of hyperplane arrangements, with an emphasis on computations
and examples. We give an introduction to many of the essential tools 
used in the area, such as Koszul and Lie algebra methods,
homological techniques, and the Bernstein-Gelfand-Gelfand correspondence,
all illustrated with concrete calculations. We also explore connections of
arrangements to other areas, such as De Concini-Procesi wonderful models,
the Feichtner-Yuzvinsky algebra of an atomic lattice, fatpoints and blowups 
of projective space, and plane curve singularities.   
\end{abstract}
\maketitle

\tableofcontents

\section{Introduction and algebraic preliminaries}
\vskip -.05in
There are a number of wonderful sources available on 
hyperplane arrangements, most notably Orlik-Terao's 
landmark 1992 text \cite{ot}. In the last decade alone several 
excellent surveys have appeared: Suciu's paper
on aspects of the fundamental group \cite{Su}, Yuzvinsky's paper on
Orlik-Solomon algebras and local systems cohomology 
\cite{Yuzsurvey}, and several monographs devoted to connections to
areas such as hypergeometric integrals \cite{othyper}, mathematical
physics \cite{var}, as well as proceedings from conferences at
Sapporo \cite{ft}, Northeastern \cite{csTapp} and Istanbul \cite{ESTUYV}.

The aim of this note is to provide an overview of some recent results 
and open problems, with a special emphasis on connections to computation. 
The paper also gives a concrete and example driven introduction
for non-specialists, but there is enough breadth here that
even experts should find something new. There are few proofs,
but rather pointers to original source material. We also explore connections of
arrangements to other areas, such as De Concini-Procesi wonderful models,
the Feichtner-Yuzvinsky algebra of an atomic lattice, the Orlik-Terao
algebra and blowups, and plane curve singularities. 
All computations in this survey can be performed using Macaulay2 \cite{GS}, 
available at: {\tt http://www.math.uiuc.edu/Macaulay2/}, and the 
arrangements package by Denham and Smith~\cite{DSmith}. 

Let $V=\K^{\ell}$, and let $S$ be the symmetric algebra on $V^*$:
$ S = \bigoplus_{i \in \Z}S_i$ 
is a $\Z$-graded ring, which means that if $s_i \in S_i$ and $s_j \in S_j$,
then $s_i \cdot s_j \in S_{i+j}$. A graded $S$-module $M$ is defined in 
similar fashion. Of special interest
is the case where $S_0$ is a field $\K$, so that each $M_i$ is a $\K$--vector space.
The free $S$ module with generator in degree $i$ is written $S(-i)$, and in 
general $M(i)_j = M_{i+j}$. 
\begin{defn} The  Hilbert function $HF(M,i) = \dim_{\K}M_i.$ \end{defn}
\begin{defn} The Hilbert series $HS(M,i) = \sum_{\Z}\dim_{\K}M_it^i.$ \end{defn}
\begin{exm}\label{exm:first}
$S=\K[x,y]$, $M=S/\langle x^2,xy \rangle$. Then
\begin{center}
\begin{supertabular}{|c|c|c|} 
\hline $i$ & $M_i$ &  $M(-2)_i$  \\ 
\hline $0$ & $1$ & $0$\\
\hline $1$ & $x,y$ &$0$\\ 
\hline $2$ & $y^2$ &$1$ \\ 
\hline $3$ & $y^3$ & $x,y$ \\ 
\hline $4$ & $y^4$ & $y^2$ \\
\hline $n$ & $y^n$ & $y^{n-2}$ \\
\hline
\end{supertabular}
\end{center}

\noindent The respective Hilbert series are
\[
HS(M,i)=\frac{1-2t^2+t^3}{(1-t)^2} \mbox{ and } HS(M(-2),i)=\frac{t^2(1-2t^2+t^3)}{(1-t)^2}
\]
An induction shows that $HS(S(-i),t)=t^i/(1-t)^{\ell}$; this makes
it easy to compute the Hilbert series of an arbitrary graded module from
a free resolution. For $S/\langle x^2,xy\rangle$, a minimal free resolution is
\[
0 \longrightarrow S(-3) \xrightarrow{\left[ \!
\begin{array}{c}
y \\
-x 
\end{array}\! \right]} S(-2)^2
\xrightarrow{\left[ \!\begin{array}{cc}
x^2& xy
\end{array}\! \right]}
 S \longrightarrow S/I \longrightarrow 0.
\]
\begin{center}
The map $[x^2,xy]$ sends 
$\begin{array}{c}
\mbox{    }e_1 \mapsto x^2\\
\mbox{    }e_2 \mapsto xy,
\end{array}$
\end{center}
so in order to have a map of graded modules, the basis elements of the
source must have degree two, explaining the shifts in the free resolution.
Taking the alternating sum of the Hilbert series yields \[
HS(M,i)=\frac{t^3-2t^2+1}{(1-t)^2}
\]
which agrees with the previous computation.
\end{exm}
\begin{exm}\label{exm:twistedcubic}
The $2\times2$ minors of $\left[ \!
\begin{array}{ccc}
x&y&z\\
y&z&w
\end{array}\! \right]$ define the twisted cubic $I \subseteq S=\K[x,y,z,w]$.
\begin{small}
\[
0 \longrightarrow S(-3)^2 \xrightarrow{\left[ \!
\begin{array}{cc}
-z & w\\
y & -z \\
-x & y
\end{array}\! \right]} S(-2)^3
\xrightarrow{\left[ \!\begin{array}{ccc}
y^2\!-\!xz& yz\!-\!xw& z^2\!-\!yw
\end{array}\! \right]}
 S \longrightarrow S/I 
\]
\end{small}
The numerical information in a free resolution may be compactly displayed 
as a {\em betti table}:
\[
b_{ij} = \dim_{\K}\Tor_i^S(M,\K)_{i+j}.
\]
$$
\vbox{\offinterlineskip 
\halign{\strut\hfil# \ \vrule\quad&# \ &# \ &# \ &# \ &# \ &# \
&# \ &# \ &# \ &# \ &# \ &# \ &# \
\cr
total&1&3&2\cr
\noalign {\hrule}
0&1 &--&--&\cr
1&--&3 &2 &\cr
\noalign{\bigskip}
\noalign{\smallskip}
}}
$$
\vskip -.1in
\noindent In particular, the indexing begins at position $(0,0)$ and is read over and down. 
So for the twisted cubic, $b_{21}(S/I) = \dim_{\K}\Tor_2^S(S/I,\K)_3 = 2.$
\end{exm}
\noindent We now give a quick review of arrangements. 
Let $\A=\{H_1,\dots ,H_n\}$ be an arrangement 
of complex hyperplanes in $\C^{\ell}$. We assume $\A$ is 
{\em central} and {\em essential}: the $\ell_i$ with $H_i = V(\ell_i)$ are
homogeneous, and the common zero locus 
$V(\ell_1,\ldots, \ell_n) = 0 \in \C^{\ell}$.
The central condition means that $\A$ also defines an 
arrangement in $\P^{{\ell}-1}$.
The main combinatorial object associated to $\A$ is the 
intersection lattice $L_{\mathcal A}$, which consists
of the intersections of elements of $\A$, ordered by reverse
inclusion. $\C^n$ is the lattice element $\hat{0}$ and the rank 
one elements of $L_{\mathcal A}$ are the hyperplanes.
\begin{defn}
The M\"{o}bius function $\mu$ : $L_{\mathcal A} \longrightarrow \Z$ is defined
 by $$\begin{array}{*{3}c}
\mu(\hat{0}) & = & 1\\
\mu(t) & = & -\!\!\sum\limits_{s < t}\mu(s) \mbox{, if } \hat{0}< t
\end{array}$$
\end{defn}
\noindent The Poincar\'e and characteristic polynomials of $\A$ are defined as
\[
\pi(\A,t) = \!\!\sum\limits_{x \in L({\mathcal A})}\mu(x) \cdot (-t)^{\text{rank}(x)}, \mbox{ and } \chi(\A,t) = t^{\text{rk}(\A)}\pi(\A,\frac{-1}{t})
\]
\begin{exm}\label{exm:nonfano}
The $A_3$ arrangement is $\bigcup_{1\le i < j \le 4}V(x_i-x_j) \subseteq \C^4$. 
Projecting along $(1,1,1,1)$ gives a central arrangement in $\C^3$, hence a
configuration of lines in $\P^2$. This configuration corresponds to the figure below, but with the line at infinity (which bounds the figure) omitted.  
\begin{center}
\epsfig{file=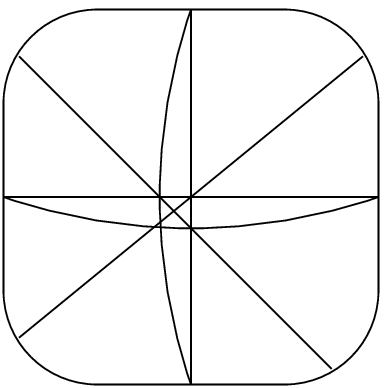,height=1.3in,width=1.3in}
\end{center}
For the $7$ rank two elements of $L(A_3)$, the four corresponding to 
triple points have $\mu = 2$, and the three normal crossings have 
$\mu =1$.  Thus, $\pi(A_3,t) = 1+6t+11t^2+6t^3$. Adding the bounding line
gives the non-Fano arrangement NF, with $\pi(NF, t)=1+7t+15t^2+9t^3$.
\end{exm}
\noindent In \cite{OS}, Orlik and Solomon showed that the cohomology ring of the
complement $\XA = \C^n\setminus \bigcup_{i=1}^d H_i$ has presentation
$H^{*}(\XA,\Z) =\bigwedge (\Z^n)/I$, with generators
$e_1, \dots , e_n$ in degree $1$ and 
\[
I = \langle \sum_{q}(-1)^{q-1}e_{i_1} \cdots \widehat{e_{i_q}}\cdots e_{i_r} \mid \codim H_{i_1} \cap \cdots \cap H_{i_r} < r \rangle.
\]
\noindent For additional background on arrangements, see \cite{ot}.
\section{$D(\A)$ and freeness}
Let $\A = \bigcup\limits_{i=1}^n H_i \subseteq V=\C^{\ell}$ be a central arrangement. For each $i$, 
fix $V(l_i)=H_i \in \A$, and define $Q_{\A} = \prod_{i=1}^n l_i \in S = \C[x_1,\ldots, x_{\ell}]$.
\begin{defn}\label{DA}
The module of $\A$-derivations (or Terao module) is the submodule of $Der_{\C}(S)$
consisting of vector fields tangent to $\A$:
\[
D({\mathcal A}) = \{ \theta \in Der_{\C}(S) | \theta(l_i) \in \langle l_i \rangle \mbox{ for all }l_i\mbox{ with }V(l_i) \in \A \}.
\]
\end{defn}
An arrangement is free when $D(\A)$ is a free $S$--module. In this case, the degrees
of the generators of $D(\A)$ are called the exponents of the arrangement. Note that
 $D(\A)$ is always nonzero, since the Euler derivation $\theta_E = \sum_{i=1}^{\ell} x_i \partial/\partial x_i \in D(\A)$. It is easy to show that 
\[
D(\A)\simeq S\cdot\theta_E \oplus syz(J_{\A}),
\]
where $J_{\A}$ is the Jacobian ideal of $Q_{\A}$, and $syz$ denotes the module of syzygies on
$J_{\A}$: polynomial relations on the generators of $J_{\A}$.
\begin{thm}[Saito \cite{slog}]
$\mathcal{A}$ is free iff there exist ${\ell}$ elements
\[
\theta_i = \sum\limits_{j=1}^{\ell} f_{ij}\frac{\partial}{\partial x_j} \in D(\mathcal{A}),
\]
such that $\det([f_{ij}]) = c \cdot Q_{\A}$, for some $c \ne 0$.
\end{thm}
\begin{exm}\label{exm:nonfanoII}
For Example~\ref{exm:nonfano}, a computation shows that
\[
\begin{array}{ccc}
D(A_3) & \simeq &S(-1)\oplus S(-2)\oplus S(-3)\\
D(NF) & \simeq &S(-1)\oplus S(-3)\oplus S(-3)
\end{array}
\]
Interestingly, 
the respective Poincar\'e polynomials factor, as
\[
\pi(A_3,t)=(1+t)(1+2t)(1+3t),\mbox{ and }\pi(NF, t)=(1+t)(1+3t)^2.
\]
This suggests the possibility of a connection between the exponents of a 
free arrangement and the Poincar\'e polynomial. 
\end{exm}
 
\noindent A landmark result in arrangements is:
\begin{thm}[Terao \cite{t}]
If $D(\A) \simeq \bigoplus\limits_{i=1}^{\ell} S(-a_i)$, then 
\[
\pi(\A,t)= \prod(1+a_it) = \sum \dim_{\C}H^i(\C^{\ell}\setminus \A)t^i.
\]
\end{thm}
\begin{exm}\label{exm:stanley}[Stanley]
For $\A$ below, $\pi(\A,t)=(1+t)(1+3t)^2$.
\begin{center}
\epsfig{file=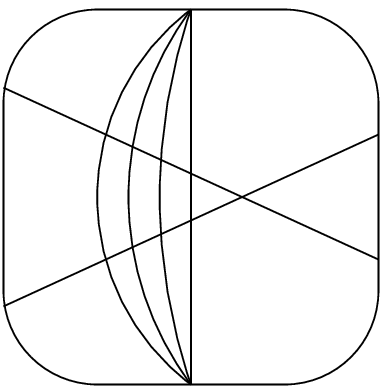,height=1.3in,width=1.3in}
\end{center}
A computation shows that $\A$ is not free, so factorization
of $\pi(\A,t)$ is a necessary but not sufficient for freeness of $\A$.
\end{exm}
A famous open conjecture in the field of arrangements is:
\begin{conj}[Terao]
If $char(\K)=0$, then 
freeness of $D(\A)$ depends only on $L_{\A}$.
\end{conj}
\begin{exm}\label{exm:zieglerAB}[Ziegler's pair \cite{z2}]
Let $\A$ be an arrangement of $9$ lines in $\P^2$, as below.
\begin{center}
\epsfig{file=yuz1.eps,height=1.5in,width=1.5in}
\end{center}
Then $D(\A)$ depends on nonlinear geometry: if the six triple points lie on a 
smooth conic, we compute:
\[
\xymatrix{0 \ar[r]& S(-7) \oplus S(-8)  \ar[r]&
S(-5) \oplus S^3(-6)  \ar[r]& syz(J_{\A})  \ar[r]& 0},
\]
while if six triple points are not on a smooth conic, 
the resolution is:
\vskip .06in
\noindent$\xymatrix{0 \ar[r]& S^4(-7)  \ar[r]&
 S^6(-6) \ar[r]& syz(J_{\A})  \ar[r]& 0}.$
\end{exm}

A version of Terao's theorem applies to any arrangement:
\begin{defn}\label{Dp}
$\Dp \subseteq \Lambda^p(Der_{\K}(S))$ consists of $\theta$ 
such that 
\[
\theta(l_i,f_2,\ldots, f_p) \in \langle l_i \rangle, \forall \mbox{ }V(l_i) \in \A, f_i \in S.
\]
\end{defn}

\begin{thm}[Solomon-Terao, \cite{solt}]\label{SolT}
\[
\chi(\A,t) = (-1)^{\ell} \lim_{x \rightarrow 1} \sum_{p \ge 0}HS(\Dp;x)(t(x-1)-1)^p.
\]
\end{thm}
\begin{prob} Relate the modules $\Dp$, for $p\ge2$, to $L_{\A}$. \end{prob}
A closed subarrangement $\hat \A \subseteq A$ is a subarrangement such that
$\hat \A = \A_X$ for some flat $X$. The best result relating 
$D(\A)$ to $L_{\A}$ is:
\begin{thm}[Terao, \cite{tunpub}] If $\hat \A \subset \A$ is a closed subarrangement, then $\pdim D(\A) \ge \pdim D(\hat \A)$.
\end{thm}
\begin{prob} Find bounds on $\pdim D(\A)$ depending on $L_{\A}$ \end{prob}
A particularly interesting class of arrangements are graphic arrangements, 
which are subarrangements of $A_n$. Given a simple (no loops or multiple edges)
graph $G$, with $\ell$ vertices and edge set $\mathsf{E}$, we define
$\A_{G}=\{z_i-z_j=0\mid (i,j)\in \mathsf{E} \subseteq \C^{\ell}\}$
\begin{thm}[Stanley \cite{stan}] $\A_{G}$ is supersolvable iff $G$ is chordal.\end{thm}
\begin{thm}[Kung-Schenck \cite{ks}] If $\A_{G}$ has an induced $k$-cycle, then $\pdim D(\A_{G})\!\ge\! k\!-\!3$. \end{thm}
\begin{exm}
The largest induced cycle of $G$ below is a 6-cycle.
\begin{center}
\epsfig{file=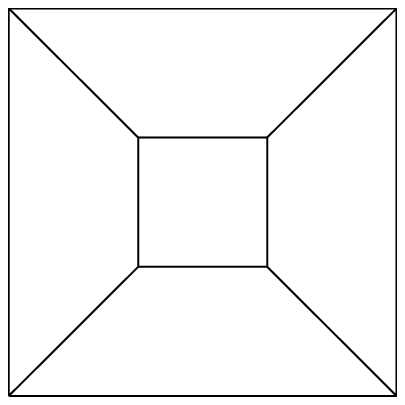,height=1in,width=1in}
\end{center}
A computation shows $\pdim(D(\A))=3$.
\end{exm}
\begin{exm}
The largest induced cycle of $G$ below is a 4-cycle.
\begin{center}
\epsfig{file=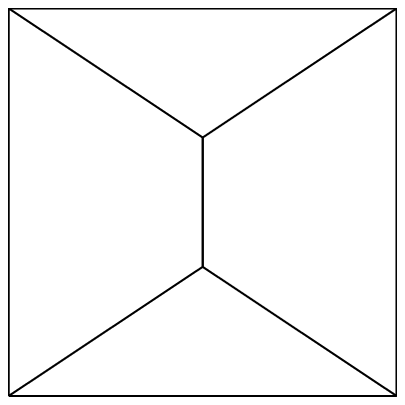,height=1in,width=1in}
\end{center}
A computation shows $\pdim(D(\A))=2$.
\end{exm}
\begin{prob} Find a formula for $\pdim D(\A_{G})$.
\end{prob}
\pagebreak
\begin{defn}\label{triple} A triple $({\mathcal A}', {\mathcal A}, {\mathcal A}'')$ of arrangements 
consists of a choice of $H \in \A$, with $\A'=\A \setminus H, \A''=\A|_H.$
\end{defn}
\noindent A main tool for proving freeness is Terao's addition-deletion theorem.
\begin{thm}[Terao \cite{t2}]\label{AD}
For a triple, any two of the following imply the third
\begin{enumerate}
\item{$D({\mathcal A})\simeq \oplus_{i=1}^n S(-b_i)$.}
\item{$D({\mathcal A}')\simeq S(-b_n +1)\oplus_{i=1}^{n-1} S(-b_i).$}
\item{$D({\mathcal A}'')\simeq \oplus_{i=1}^{n-1} S/L(-b_i).$}
\end{enumerate}
\end{thm}
\begin{exm}
In Example~\ref{exm:nonfano}, the $A_3$ arrangement is free with exponents $\{1,2,3\}$.
Let $H$ be the line at infinity, which meets $A_3$ in four points. Then $D({\mathcal A}'')$ 
is free, with exponents $\{1,2\}$, so the non-Fano arrangement is free with exponents $\{1,3,3\}$, which 
agrees with our earlier computation. Example~4.59 of \cite{ot} gives a 
free arrangement for which the addition-deletion theorem does not apply.
\end{exm}
As a corollary of Theorem~\ref{AD}, Terao showed that supersolvable
arrangements are free.
\begin{defn}\label{SSolve}
An element $X$ of a lattice is modular if for all $Y \in L$ and all $Z<Y$,
$Z \vee (X \wedge Y) = (Z \vee X) \wedge Y$. 
A central arrangement $\A$ is supersolvable if there exists a maximal chain 
$\hat{0} = X_0 < X_1 < \cdots < X_n = \hat{1}$ of modular elements in $L(\A)$.
\end{defn}
For line configurations
in $\P^2$, the supersolvability condition simply means there is a singular point
$p \in \A$ such that every other singularity of $\A$ lies on a line of $\A$ which
passes through $p$. For example, the $A_3$ arrangement is supersolvable, since 
any triple point is such a singularity. For arrangements in $\P^2$, there
is a beautiful characterization of freeness involving multiarrangements.
\begin{defn}
A multiarrangement $(\mathcal{A},{\bf m})$ consists of an arrangement $\A$,
along with a multiplicity $m_i \in \N$ for each $H \in \A$.
\[
D(\A,{\bf m}) = \{\theta \mid \theta(l_i) \in \langle l_i^{m_i}\rangle \}.
\]
\end{defn}
%
  
\begin{thm}
$\A \subseteq \mathbb{P}^2$ is free if and only if 
\begin{enumerate}
\item{$\pi(\A,t) = (1+t)(1+at)(1+bt)$ and}
\item{$D(\A|_H,{\bf m}) \simeq S/L(-a)\oplus S/L(-b)$,} 
\end{enumerate}
where $(2)$ holds for all $H=V(L) \in \A$, with ${\bf m}(H_i)\!=\!\mu_{\A}(H\cap H_i).$
\end{thm}
The necessity of these conditions was shown by Ziegler in \cite{z},
and sufficiency was proved by Yoshinaga in \cite{y2}. 
In \cite{y1}, Yoshinaga gives a generalization to higher dimensions.
\section{Multiarrangements}
\noindent The exponents of free multiarrangements are not combinatorial:
\begin{exm}[Ziegler, \cite{z}]
Consider the two multiarrangements in $\P^1$, with underlying arrangements
defined by 

\begin{center}
$\begin{array}{cccc}
\A_1 &= &V(x\cdot y \cdot (x+y)\cdot (x-y)) \\
\A_2 &= &V(x\cdot y \cdot (x+y)\cdot (x-ay)),
\end{array}$
\end{center}
with $a \ne 1$. To compute $D(\A_1,(1,1,3,3))$, we must find all
\[
\theta = f_1(x,y)\partial/\partial x + f_2\partial/\partial y
\]
such that 
\[
\begin{array}{c}
\theta(x) \in \langle x \rangle, \mbox{ } \theta(x+y) \in \langle x+y \rangle^3\\
\theta(y) \in \langle y \rangle, \mbox{ }  \theta(x-y) \in \langle x-y \rangle^3
\end{array}
\]
Thus, $D(\A_1,(1,1,3,3))$ is the kernel of the matrix 
\begin{center}
$
\left[ \!
\begin{array}{cccccc}
1 &0 & x& 0  & 0 &                              0     \\
0&1& 0&  y&0 &                                  0 \\  
1& 1& 0&  0&  (x+y)^3& 0\\
1&-1&0&  0&  0&                              (x-y)^3 
\end{array}\! \right].$
\end{center}
Computations show that $D(\A_1,(1,1,3,3))$ has exponents $\{3,5\}$,
and $D(\A_2,(1,1,3,3))$ has exponents $\{4,4\}$. 
\end{exm}
\noindent There is an analog of Theorem~\ref{SolT} for multiarrangements.
\begin{defn}\label{Dp2}
$D^p(\A,{\bf m}) \subseteq \Lambda^p(Der_{\K}(S))$ consists of $\theta$ 
such that 
\[
\theta(l_i,f_2,\ldots, f_p) \in \langle l_i \rangle^{m(\l_i)}, \forall \mbox{ }V(l_i) \in \A, f_i \in S.
\]
\end{defn}
\begin{thm}[Abe-Terao-Wakefield \cite{atw2}]
Define 

$\begin{array}{ccc}
\Psi(\A,{\bf m},t,x) &=&\sum\limits_{p=0}^{\ell} HS(D^p(\A,{\bf m}),x)(t(x-1)-1)^p \\
\chi((\A,{\bf m}),t) &=&(-1)^{\ell}\lim_{x \rightarrow 1}\Psi(\A,{\bf m},t,1).
\end{array}$

If $D^1(\A,{\bf m}) \simeq \oplus S(-d_i),$ then $\chi((\A,{\bf m}),t)= \prod\limits_{i=1}^{\ell}(1+d_it).$
\end{thm}
In \cite{atw}, Abe-Terao-Wakefield prove an addition-deletion theorem for multiarrangements by introducing {\em Euler multiplicity} for the restriction. 
It follows from the Hilbert-Burch 
theorem that any $(\A,{\bf m}) \subseteq \P^1$ is free, which leads to the question
of whether there exist other arrangements which are free for any ${\bf m}$. In \cite{aty},
Abe-Terao-Yoshinaga prove that any such arrangement is a product of one- and two--dimensional
arrangements. Nevertheless, several natural questions arise:
\begin{prob} Characterize the projective dimension of $D(\A,{\bf m})$.
\end{prob}
\begin{prob} Define supersolvability for multiarrangements.
\end{prob}
\section{Arrangements of plane curves}
For a collection of hypersurfaces 
\[
\CL=\bigcup_i V(f_i) \subseteq \P^n,
\]
the
module of derivations $D(\CL)$ is obtained by substituting $f_i$ for $l_i$ in Definition~\ref{DA}. 
It is not hard to prove that
Saito's criterion still applies. Are there other freeness theorems? 
\begin{exm}
For the arrangement $\CL \subseteq \P^2$ depicted below
\begin{center}
\epsfig{file=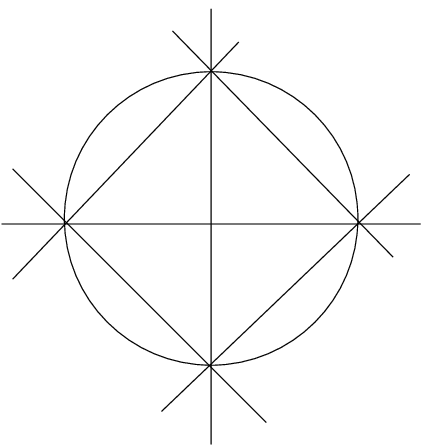,height=1.5in,width=1.5in}
\end{center}
we compute that $D(\CL) \simeq S(-1) \oplus S(-2) \oplus S(-5)$.
\end{exm}
This example can be explained by an addition-deletion
theorem \cite{ST1}, but there is subtle behavior related to singular points. For the remainder
of this section, $C=\cup_i V(f_i)\subseteq \C^2$ is reduced plane curve, and if $p\in C$ is a singular
point, translate so $p=(0,0)$. 
\begin{defn}\label{qhomogSing}
A plane curve singularity is quasihomogeneous if and only if
there exists a holomorphic change of variables so that 
$f(x,y) = \sum c_{ij}x^{i} y^{j}$ is weighted homogeneous: there exists 
$\alpha, \beta \in \Q$ such that $\sum c_{ij}x^{i \cdot \alpha} y^{j \cdot \beta}$ is homogeneous.
\end{defn}
  
\begin{defn}
The Milnor number at $(0,0)$ is 
$$
\mu_{(0,0)}(C) = \dim_\mathbb{C}\mathbb{C}\{x,y\}/\langle\frac{\partial f}{\partial x}\mbox{, }\frac{\partial f}{\partial y}\rangle.
$$
The Tjurina number at $(0,0)$ is 
$$
\tau_{(0,0)}(C) = \dim_{\C} \C \{x,y\}/\langle\frac{\partial f}{\partial
x}\mbox{, }\frac{\partial f}{\partial y}\mbox{, }f\rangle.
$$
\end{defn}
For a projective plane curve $V(Q)\subseteq \P^2$, it is easy 
to see that the degree of $Jac(Q) = \sum_{p \in sing(V(Q))}\tau_p$.
\begin{exm}
Let $\CL$ be as below:
\begin{center}
\epsfig{file=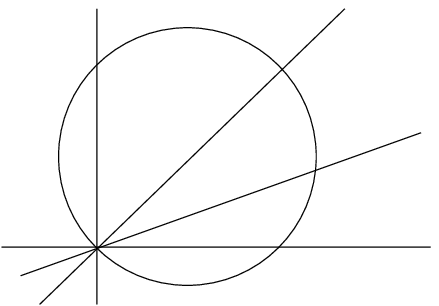,height=1.1in,width=1.6in}
\end{center}
If $p$ is an ordinary singularity with $k$ distinct branches, 
then $\mu_p(C)=(k-1)^2$, so the sum of the
Milnor numbers is $20$. However, a computation shows that
$\deg(J_{\CL}) = 19$. All singularities are ordinary, but the 
singularity at the origin is not quasihomogeneous.
\end{exm}

\begin{thm}[Saito \cite{s}]
If $C=V(f)$ has an isolated singularity at the origin, then $f \in Jac(f)$ iff $f$ is quasihomogeneous.
\end{thm}
For arrangements of lines and conics such that every singular point is 
quasihomogeneous, \cite{ST1} proves an addition/deletion theorem; 
\cite{STY} generalizes the result to curves of higher genus.
\begin{exm}
Let $\CL$ be as below:
\begin{center}
\epsfig{file=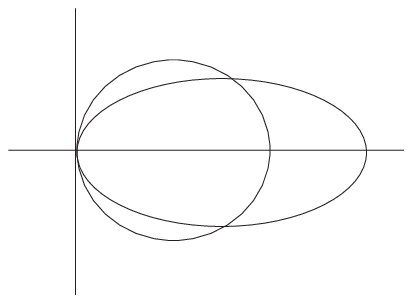,height=1.5in,width=1.5in}
\end{center}
$D(\CL)$ has exponents $\{1,2,3\}$, which can be shown using the aforementioned addition-deletion
theorem. Change $\CL$ to $\CL'$ via: 
\[
y=0 \longrightarrow x-13y=0.
\]
A computation shows that
$D(\CL')$ is not free. Thus, for line-conic arrangements, freeness is not combinatorial.
\end{exm}
\begin{prob} Define supersolvability for hypersurface arrangements.
\end{prob}
\begin{prob} Give combinatorial bounds on $\pdim D(\CL)$.
\end{prob}
\begin{prob} Analyze associated primes and $\Ext$ modules of $D(\CL)$.
\end{prob}

\section{The Orlik--Terao algebra and blowups}
The Orlik--Terao algebra is a symmetric analog of
the Orlik-Solomon algebra. While the Orlik-Solomon algebra
records the existence of dependencies among sets of hyperplanes, the Orlik-Terao
algebra records the actual dependencies. If
$\codim \cap_{j=1}^m H_{i_j} < m$, then there exist $c_{i_j}$ with
\[
\sum\limits_{j=1}^m c_{i_j}\cdot l_{i_j} = 0 \mbox{ a dependency}.
\]
\begin{defn}
The Orlik-Terao ideal 
\[
I_{\A} = \langle \sum_{j=1}^m c_{i_j} (y_{i_1}\cdots \hat y_{i_{j}} \cdots y_{i_m}) \mid \mbox{ over all dependencies}\rangle
\]
\end{defn}
The Orlik-Terao algebra is $C(\A) = \K[x_1,\ldots,x_n]/I_{\A}.$
\begin{exm}\label{otexm1}
$\A = V(x_1 \cdot x_2\cdot x_3 \cdot (x_1+x_2+x_3))$, the only dependency is
$l_1\!+\!l_2\!+\!l_3\!-\!l_4 = 0$, so $I_{\A} =\langle y_2y_3y_4+y_1y_3y_4+y_1y_2y_4-y_1y_2y_3 \rangle.$
\end{exm}
In \cite{ot1}, Orlik and Terao answer a question of Aomoto by considering
the Artinian quotient $AOT$ of $C(\A)$ by $\langle x_1^2,\ldots,x_n^2 \rangle$. They prove:
\begin{thm}[Orlik-Terao \cite{ot1}] $HS(AOT, t) = \pi(\A,t)$.\end{thm}
\begin{thm}[Terao \cite{ter}]\label{TeraoCA}
\[
HS(C(\A),t) = \pi\Big(\A, \frac{t}{1-t}\Big).
\]
\end{thm}
It is not hard to show that 
\[0 \rightarrow I_{\A} \rightarrow \K[x_1,\ldots,x_n]\stackrel{\phi}{\rightarrow} \K\Bigg[\frac{1}{l_1},\ldots, \frac{1}{l_n}\Bigg] \rightarrow 0\]
is exact, so $V(I_{\A})\subseteq \P^{n-1}$ is irreducible and rational. In
any situation where weights of dependencies play a role, the Orlik-Terao
algebra is the natural candidate to study. One such situation involves 2-formality:
\begin{defn}
$\A$ is 2-formal if all dependencies are generated by
dependencies among three hyperplanes.
\end{defn}
\begin{thm}[Falk-Randell \cite{fr}]If $\A$ is $K(\pi,1)$, $\A$ is 2-formal.\end{thm}
\begin{thm}[Yuzvinsky \cite{yu1}] If $\A$ is free, $\A$ is 2-formal.\end{thm}
One reason that formality is interesting is that it is not a combinatorial
invariant: in Example~\ref{exm:zieglerAB}, the arrangement 
for which the six triple points lie on a smooth conic is not 
2-formal, and the arrangement for which the points do not lie on a 
smooth conic is 2-formal. 
\begin{thm}[\cite{ST1}]\label{ST2formal} $\A$ is 2-formal iff $\codim(I_{\A})_2=n-\ell.$
\end{thm}
In \cite{bt}, Brandt and Terao generalized the
notion of 2-formality to $k-$formality: 
$\A$ is k-formal if all dependencies are generated by
dependencies among $k+1$ or fewer hyperplanes.
Brandt and Terao prove that every free arrangement is $k-$formal.
\begin{prob}
Find an analog of Theorem~\ref{ST2formal} for $k-$formality.
\end{prob}
\begin{exm}
The configuration of Example~\ref{otexm1} consists
of four generic lines:
\vskip -.1in
\begin{figure}[ht]
\subfigure{%
\begin{minipage}[t]{0.5\textwidth}
\setlength{\unitlength}{17pt}
\begin{picture}(4.5,4.5)(-2.5,-0.5)
\put(0,0){\line(1,1){4}}
\put(-1,2){\line(1,0){6}}
\put(.5,4){\line(1,-1){4}}
\put(-0.5,0.7){\line(3,1){5.5}}
\end{picture}
\end{minipage}
}
\end{figure}
\vskip -.1in
\noindent The Orlik-Terao ideal defines a cubic surface in $\P^3$, and
a computation shows that $V(I_{\A})$ has four singular points.
\end{exm}
This can be interpreted in terms of a rational map. Let $\alpha_i = Q_{\A}/l_i$,
and define $\phi_{\A}=[\alpha_1,\ldots,\alpha_n]$.  
\[
\P^{\ell-1}\stackrel{\phi_{\A}}{\longrightarrow}\P^{n-1},
\]
Restrict to the case $\A \subseteq \P^2$, and let 
$\xA \stackrel{\pi}{\longrightarrow} \P^2$ denote the blowup of 
$\P^2$ at the singular points of $\A$, with $E_0$ denoting the 
pullback to $\xA$ of the class of a line on $\P^2$, and $E_i$ 
the exceptional divisors over singular points of $\A$. Let 
\[
D_{\A} = (n-1)E_0 - \!\!\sum\limits_{p_i \in L_2(\A)} \mu(p_i)E_i.
\]
Utilizing results of Proudfoot-Speyer \cite{ps} showing that $C(\A)$ is 
Cohen-Macaulay and the Riemann-Roch theorem, \cite{syzRes} 
shows that the map $\phi_{\A}$ is determined by the global sections 
of $D_{\A}$, and that $\phi_{\A}$
\begin{enumerate}
\item is an isomorphism on $\pi^*(\P^2 \setminus \A$)
\item contracts the lines of $\A$ to points
\item blows up the singularities of $\A$.
\end{enumerate}

\pagebreak

\begin{defn}\label{CMregularity}
A graded $S$-module $N$ has Castelnuovo-Mumford 
regularity $d$ if $\Ext^j(N,S)_n=0$ for all $j$ and all $n\le -d-j-1$. 
\end{defn}
In terms of the betti table, the regularity of $N$ is the label of the 
last non-zero row, so in Example~\ref{exm:twistedcubic}, $S/I$ has
Castelnuovo-Mumford regularity one. The regularity of $C(\A)$
is determined in \cite{syzRes}:
\begin{thm}[\cite{syzRes}]\label{CMregofCA}
For $\A \subseteq \P^{\ell-1}$, $C(\A)$ is $\ell-1$--regular.
\end{thm}
To see this, note that since $C(\A)$ is Cohen-Macaulay, 
quotienting $C(\A)$ by $\ell$ generic linear forms yields an Artinian 
ring whose Hilbert series is the numerator of the Hilbert 
series of $C(\A)$. The regularity of an Artinian module is 
equal to the length of the module, so the result follows from 
Theorem~\ref{TeraoCA}.

A main motivation for studying $C(\A)$ is a surprising connection 
to nets and resonance varieties, which are the subject of 
\S~\ref{resonanceSection}. First, the definition of a net:

\begin{defn}~\label{Net}
Let $3 \le k \in \Z$. A $k$-net in $\P^2$ is a partition of the lines of 
an arrangement $\A$ 
into $k$ subsets $\A_i$, together with a choice of points $Z \subseteq \A$, such that:
\begin{enumerate}
\item{for every $i \ne j$ and every $L \in\A_i, \ L'\in\A_j$, $L \cap L'\in Z$.}
\item{$\forall$ $p \in Z$ and every $i \in \{1,\ldots,k\}$, there exists a unique $L\in A_i$ with $Z\in L$.}
\end{enumerate}
\end{defn}

In \cite{LY}, Libgober and Yuzvinsky show that nets are related to 
the first resonance variety $R^1(\A)$. The definition of a net 
forces each subset $\A_i$ to have the same 
cardinality, and if $m=|\A_i|$, the net is called a
$(k,m)$-net. Using work of \cite{LY} and \cite{FalkY}, it is 
shown in \cite{syzRes} that
\begin{thm}\label{ENcpx} Existence of a $(k,m)$ net 
implies that there is a decomposition $D_{\A} = A+B$ with $h^0(A)=2$
and $h^0(B)=km - {m+1 \choose 2}$.
\end{thm}
\begin{defn}\label{1generic}
A matrix of linear forms is $1$-generic if it has no zero entry, and
cannot be transformed by row and column operations to have a zero entry.
\end{defn}
In \cite{EisenbudFactor}, Eisenbud shows that if a divisor $D$ on a
smooth curve $X$ factors as $D \simeq A+B$, with $A$ having 
$m$--sections and $B$ having $n$--sections, then the ideal of the
image of $X$ under the map defined by the global sections
of $D$ will contain the $2 \times 2$ minors of a $1$-generic matrix.
Using this result and Theorem~\ref{ENcpx}, it can be shown that
$I_{\A}$ contains the ideal $I_2(M)$ of $2 \times 2$ minors of a 
$1$-generic \begin{small}$2 \times \Big( km - {m+1 \choose 2} \Big)$\end{small} matrix $M$.
So if \begin{small}$G=S(-1)^{km - {m+1 \choose 2}}$\end{small}, the 
Eagon-Northcott complex \cite{E}
\[
\cdots \rightarrow S_2(S^2)^* \otimes \Lambda^4 G \rightarrow (S^2)^* 
\otimes \Lambda^3 G \rightarrow \Lambda^2 G \rightarrow \Lambda^2 S^2 \rightarrow S/I_2(M) \rightarrow 0
\]
is a subcomplex of resolution of $S/I_{\A}$. The geometric 
content of Theorem~\ref{ENcpx} is that it implies $V(I_{\A})$ 
lies on a {\em scroll} \cite{E}. 

\begin{exm}\label{ShowMeThm}
For the $A_3$ arrangement, the set of triple points $Z$ gives a 
$(3,2)$ net, where the $A_i$ correspond to normal crossing points:
$A_1 = 12\mid 34$, $A_2 = 13\mid 24$, $A_3 = 14\mid 23$.
\begin{figure}[ht]
\subfigure{%
\begin{minipage}[t]{0.30\textwidth}
\setlength{\unitlength}{16pt}
\begin{picture}(6,7)(-1,-2.5)
\put(-0.8,0){\line(1,0){5.6}}
\put(-0.8,-0.4){\line(2,1){4.7}}
\put(-0.4,-0.8){\line(1,2){2.8}}
\put(2,-0.8){\line(0,1){5.6}}
\put(4.8,-0.4){\line(-2,1){4.7}}
\put(4.4,-0.8){\line(-1,2){2.8}}
\put(-1.7,0.3){\makebox(0,0){$L_{12}$}}
\put(-1.5,-0.6){\makebox(0,0){$L_{13}$}}
\put(-0.5,-1.3){\makebox(0,0){$L_{23}$}}
\put(2,-1.3){\makebox(0,0){$L_{34}$}}
\put(4.5,-1.3){\makebox(0,0){$L_{24}$}}
\put(5.5,-0.6){\makebox(0,0){$L_{14}$}}
\end{picture}
\end{minipage}
}
\end{figure}
\[
\mbox{Let }A = 2E_0-\!\!\sum\limits_{\{p | \mu(p) = 2\}} \!\!E_p \mbox{ and } B = 3E_0-\!\!\sum\limits_{p \in L_2(\A)} \!\!E_p.
\]
So $n - {m+1 \choose 2} = 6-3 =3$ and $I$ 
contains the $2\times 2$ minors of a $2 \times 3$ matrix, whose
resolution appears in Example~\ref{exm:twistedcubic}. 
The graded betti diagram for $\C[x_0,\ldots,x_5]/I_{A}$ is
$$
\vbox{\offinterlineskip 
\halign{\strut\hfil# \ \vrule\quad&# \ &# \ &# \ &# \ &# \ &# \
&# \ &# \ &# \ &# \ &# \ &# \ &# \ &# \
\cr
total&1&4&5&2&\cr
\noalign {\hrule}
0&1 &--&--&--&\cr
1&--&4 &2 &--&\cr
2&--&--&3 &2&\cr
\noalign{\bigskip}
\noalign{\smallskip}
}}
$$
From this, it follows that the free resolution of $S/I_{A}$ 
is a mapping cone resolution \cite{E}. The geometric meaning is that 
$X_{\A}$ is the intersection of a generic quadric hypersurface with 
the scroll. \end{exm}
\noindent Since $D_{\A}$ contracts proper transforms of lines to points,
it is not very ample. However, it follows from \cite{syzRes} that 
$D_{\A}+E_0$ is very ample, 
and gives a De-Concini-Procesi wonderful model (see next section) 
for the blowup. 
\begin{prob} Determine the graded betti numbers of $C(\A)$.
\end{prob}
 \begin{prob} Relate $R^k(\A)$ to the graded betti numbers of $C(\A)$.
\end{prob}

\section{Compactifications} 
In \cite{FM}, Fulton-MacPherson provide a compactification $F(X,n)$ for
the configuration space of $n$ marked points on an algebraic variety $X$. 
The construction is quite involved, but the combinatorial data is that of 
$A_n$. In a related vein, in \cite{DP}, De Concini-Procesi construct 
a wonderful model $X$ for a subspace complement 
$\XA = \C^{\ell} \setminus \A$: a smooth, compact $X$ such that 
$X \setminus \XA$ is a normal crossing divisor. Here
it is the combinatorics which are complex. A key object in their construction is
\[
\XA \longrightarrow \C^{\ell} \times \prod\limits_{D \in G}\P(\C^{\ell}/D),\]
where $G$ is a building set. In \cite{FK}, Feichtner-Kozlov generalize the construction of \cite{DP} to a purely lattice-theoretic setting. See \cite{F} for 
additional background on this section.
\begin{defn}\label{definition:NestSet}
For a lattice $L$, a building set $G$ is a subset of $L$, such that
for all $x \in L$, 
$\max \{G_{\le x}\} = \{x_1,\ldots,x_m\}$ satisfies $[\hat{0},x] \simeq \prod_{j=1}^m [\hat{0},x_j].$
A building set contains all irreducible $x \in L$.
\end{defn}
\begin{defn}\label{definition:NestSet}
A subset $N$ of a building set $G$ is nested if for any 
set of incomparable $\{x_1,\ldots,x_p\} \subseteq N$ with 
$p\ge 2$, $x_1 \vee x_2 \vee \cdots \vee x_p$ exists in $L$, but 
is not in $G$. 
\end{defn}
Nested sets form a simplicial complex $N(G)$, with
vertices the elements of $G$ (which are vacuously nested). 
 \begin{exm}\label{minBldgSet}
The minimal building set for $A_3$ consists of the hyperplanes themselves, the
triple intersections in $L_2$, and the element $\hat{1}$. Since $\hat{1}$
is a member of every face of $N(G)$, the nested set complex $N(G)$ is the cone over 
\begin{center}
\epsfig{file=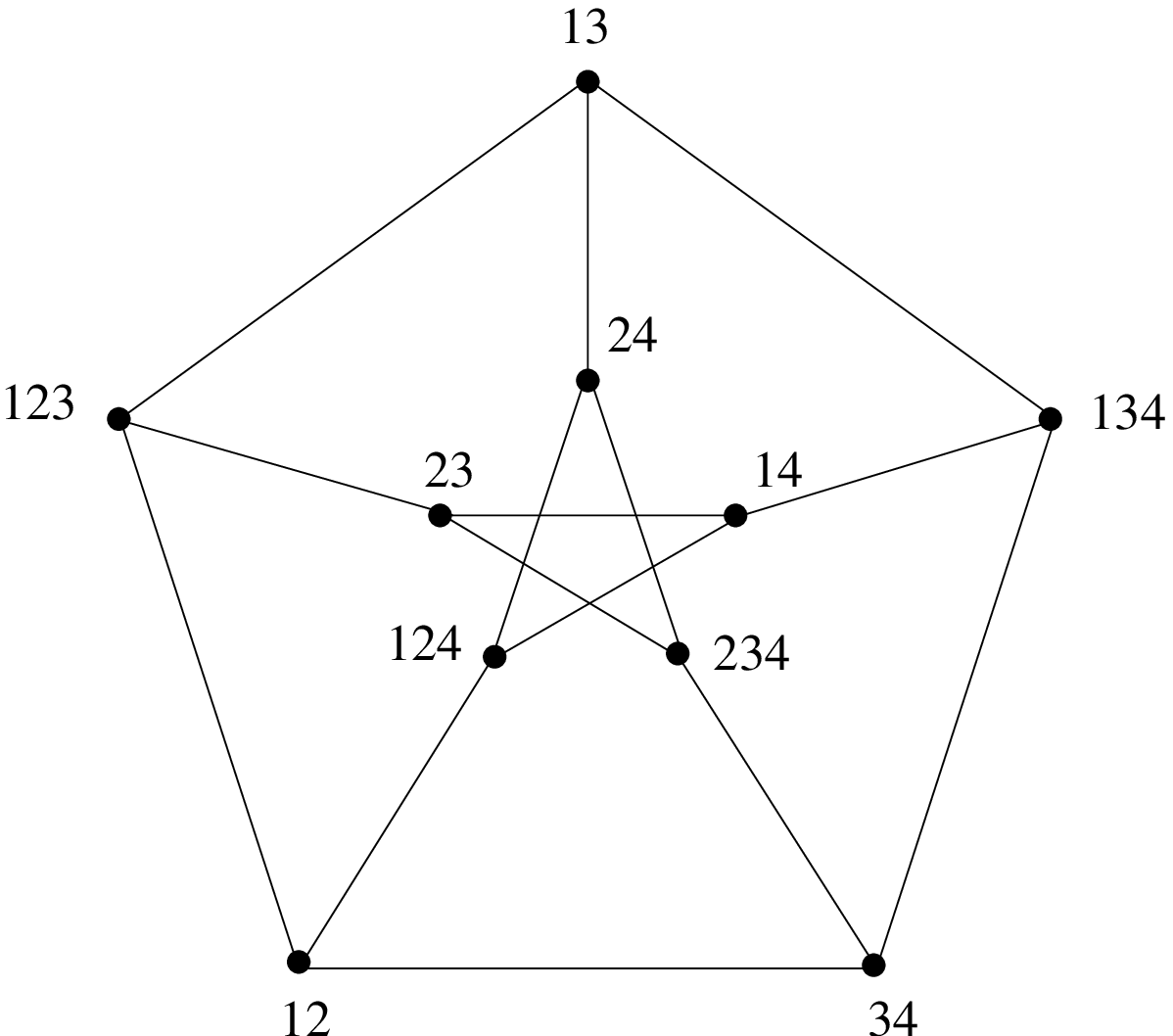,height=1.7in,width=2.0in}
\end{center}
There is an edge $\overline{(12),(123)}$ because there are no incomparable subsets with 
at least two elements, while $\overline{(12),(34)}$ is an edge because $(12)\vee(34)$ exists
in $L$ (it is a normal crossing), but is not in $G$.
\end{exm}
Suppose $L$ is an atomic lattice, and $G$ a building set in $L$. In
\cite{FY}, Feichtner and Yuzvinsky study a certain algebra associated to
the pair $L,G$:
\[
D(L,G) = \Z[x_g | g \in G]/I, \mbox{ with }x_g \mbox{ of degree }2.
\]
where $I$ is generated by 
\[
\prod\limits_{\{g_1,\ldots, g_n\} \not\in N(G)\}} x_{g_i} \mbox{ and } 
\sum\limits_{g_i \ge H \in L_1}x_{g_i}
\]
\begin{thm}[Feichtner-Yuzvinsky \cite{FY}]\label{thm:AtomicLatticealg}
If $\A$ is a hyperplane arrangement and $G$ a building set 
containing $\hat{1}$, then 
\[
D(L,G) \simeq H^*(Y_{\A,G},\Z),
\]
where $Y_{\A,G}$ is the wonderful model arising from the building set $G$.
\end{thm}
The importance of this is the relation to the
Knudson-Mumford compactification $\overline{M_{0,n}}$ of the moduli space
of $n$ marked points on $\P^1$.
\begin{thm}[De Concini-Procesi \cite{DP2}]\label{dpthm}
\[
\overline{M_{0,n}} \simeq Y_{A_{n-2}, G},
\]
where $G$ is the minimal building set for $A_{n-2}$. 
\end{thm}
A presentation for the cohomology ring 
of $\overline{M_{0,n}}$ was first described by Keel in \cite{keel}; the 
description which follows from \cite{FY} is very economic.
\begin{exm}
By Theorem~\ref{thm:AtomicLatticealg} and \cite{DP2}, 
\[
H^*(\overline{M_{0,5}}, \Z) \simeq D(L(A_3), G_{min}).
\]
The nested set complex for $A_3$ and $G_{min}$ appears in Example~\ref{minBldgSet},
so that $D(L(A_3), G_{min})$ is the quotient of a polynomial ring $S$ with eleven generators
by an ideal consisting of 6 linear forms (one form for each hyperplane) and 19 quadrics.
To see that there are 19 quadrics, note that the space of quadrics in $11$-variables 
has dimension $45$, and $N(G_{min})$ has $15+11 = 26$ edges (recall that $\hat{1}$ is
not pictured). A computation shows that
\[
D(L(A_3), G_{min})\simeq \Z[x_1,\ldots, x_5]/I,
\]
where $I$ consists of all but one quadric of $S$ (and includes all squares of variables).
This meshes with the intuitive picture: to obtain a wonderful model, simply blow up the
four triple points, so that $\overline{M_{0,5}}$ is the corresponding Del Pezzo surface $X_4$, which
has $\sum h^i(X_4,\Z)t^i = 1+5t^2+t^4$, agreeing with the computation.
\end{exm}
\begin{prob} Analyze $D(L,G)$ for other lattices.
\end{prob}
 \section{Associated Lie algebra of $\pi_1$ and LCS ranks}
Let $G$ be a finitely-generated group, with
normal subgroups, 
\[
G=G_1 \ge G_2 \ge G_3\ge\cdots,
\]
defined inductively by $G_k = [G_{k-1},G]$.
We obtain an associated Lie algebra
\[
gr(G)\otimes \Q := \bigoplus_{k=1}^{\infty} G_k/G_{k+1} \otimes \Q,
\]
with Lie bracket induced by the commutator map.
Let $\phi_k=\phi_k(G)$ denote the rank of the $k$-th quotient.
Presentations for $\pi_1(\XA)$ are given by Randell \cite{R}, Salvetti \cite{Sal}, Arvola \cite{Ar},
and Cohen-Suciu \cite{CS}. For computations, the braid monodromy presentation of \cite{CS} 
is easiest to implement. For a detailed discussion of $\pi_1(\XA)$, see 
Suciu's survey \cite{Su}.
The fundamental group is quite delicate, and in this section, we investigate
properties of $\pi_1(\XA)$ via the associated graded Lie algebra
\[
\g = gr(\pi_1(\XA))\otimes \Q
\]
The Lefschetz-type theorem of Hamm-Le \cite{HL} implies that taking 
a generic two dimensional slice gives an isomorphism on $\pi_1$.
Thus, to study $\pi_1(\XA)$, we may assume $\A \subseteq \C^2$ or $\P^2$.
As shown by Rybnikov \cite{Ry}, $\pi_1(\XA)$ is not determined by $L_{\A}$;
whereas the Orlik-Solomon algebra $H^*(\XA,\Z)$ is determined by 
$L_{\A}$.
\begin{exm}\label{A3HS}
In Example~\ref{exm:nonfano}, we saw that the Hilbert series for $A_3$ is
$1+6t+11t^2+6t^3$. A computation shows that the LCS ranks begin
\begin{center}
$\begin{array}{cccccc}
6 &4 & 10 &21 & 54 & \cdots
\end{array}$
\end{center}
For higher $k$, $\phi_k(\pi_1(A_3)) = w_k(2)+w_k(3)$, where $w_k$ is a Witt number. 
In general, we may encode the LCS ranks via
\[
\prod_{k=1}^{\infty}\frac{1}{(1-t^k)^{\phi_k}} 
\]
For $A_3$, this is
\[
\frac{1}{(1\!-\!t)^6}\frac{1}{(1\!-\!t^2)^4}\frac{1}{(1\!-\!t^3)^{10}}\frac{1}{(1\!-\!t^4)^{21}}\frac{1}{(1\!-\!t^5)^{54}}\cdots
\]
Expanding this and writing out the first few terms yields
\[
1+ 6t +25t^2+ 90t^3 +301t^4+966t^5+3025t^6+\cdots
\]
\vskip .1in
If we multiply this with
\[
\pi(A_3,-t)=  1- 6t +11t^2-6t^3,
\]
the result is $1$, and is part of a general pattern. 
\end{exm}
\begin{thm}[Kohno's LCS formula \cite{K85}]
For the arrangement $A_{n-1}$ (graphic arrangement of $K_n$)
\[
\prod_{k=1}^{\infty}(1-t^k)^{\phi_k}=\prod\limits_{i=1}^{n-1}(1-it).
\]
\end{thm}
This explains the computation of Example~\ref{A3HS}. 
We now compute the free resolution of the residue field $A/m$ as an $A$-module, where
$m = \langle E_1 \rangle$. Let
\[
b_{ij} = \dim_{\Q}Tor^{A}_i(\Q,\Q)_j 
\]
\begin{exm}\label{A3Tors}
For $A_3$, we compute $b_{ij}=0$ if $i \ne j$, and 
\[
\sum_ib_{ii}t^i = 1+ 6t +25t^2+ 90t^3 +301t^4+966t^5+3025t^6+\cdots
\]
The $b_{ii}$ are the coefficients of the formal power series in Example~\ref{A3HS}!
\end{exm} 
Kohno's result was the first of a long line of results on 
LCS formulas for certain special families of arrangements
\begin{enumerate}
\item{Braid arrangements: Kohno \cite{K85}}
\item{Fiber type arrangements: Falk--Randell \cite{FR85}}
\item{Supersolvable arrangements: Terao \cite{T}}
\item{Lower bound for $\phi_k$: Falk \cite{Fa89}}
\item{Koszul arrangements: Shelton--Yuzvinsky \cite{SY}}
\item{Hypersolvable arrangements: Jambu--Papadima \cite{JP}}
\item{Rational $K(\pi,1)$ arrangements: Papadima--Yuzvinsky \cite{PY}}
\item{MLS arrangements: Papadima--Suciu \cite{PS}}
\item{Graphic arrangements: Lima-Filho--Schenck \cite{LS}}
\item{No such formula in general: Peeva \cite{P}}
\end{enumerate}
Let $\LL(H_1(\XA,\K))$ denote the free Lie algebra on $H_1(\XA,\K).$
Dualizing the cup product gives a map 
\[
H_2(\XA,\Q) \stackrel{c}{\rightarrow} H_1(\XA,\Q)\wedge H_1(\XA,\Q)
\longrightarrow \LL(H_1(\XA,\Q)),
\]
Following Chen \cite{C}, define the holonomy Lie algebra
\[
\hA = \LL(H_1(\XA,\K))/I_{\A},
\]  
where $I_{\A}$ is the Lie ideal generated by $\im(c)$. As noted by Kohno in
\cite{K83}, taking transpose of cup product shows that the 
image of $c$ is generated by 
\[
[x_j, \sum_{i=1}^k x_i],
\]
where $x_i$ is a generator of $\LL(H_1(X,\K))$ corresponding to $H_i$, and
the set $\{H_1,\ldots, H_k \}$ is a maximal dependent set of codimension two, 
so corresponds to an element of $L_2(\A)$. The upshot is that
\[
\prod_{k=1}^{\infty} \frac{1}{(1-t^k)^{\phi_k}} = \sum_{i=0}^{\infty}\dim_{\Q} Tor_i^A(\Q,\Q)_it^i.
\]
This was first made explicit by Peeva in \cite{P}; the proof runs as follows.
First, Brieskorn \cite{Br} showed that $\XA$ is formal, in the sense of \cite{S}. 
Using Sullivan's work and an analysis of the bigrading on Hirsch extensions,
Kohno proved
\begin{thm}[Kohno] $\phi_k(\g) = \phi_k(\hA)$.\end{thm}
Thus
\begin{enumerate}
\item{$\prod_{k=1}^{\infty} \frac{1}{(1-t^k)^{\phi_k}} = HS(U(\hA,t))$, 
which follows from Kohno's work and Poincar\'e-Birkhoff-Witt.}
\item{Shelton-Yuzvinsky show in \cite{SY} that $U(\hA)=\overline{A}^{!}$ is the quadratic dual of the quadratic Orlik-Solomon algebra.}
\item{Results of Priddy-L\"ofwall show that the quadratic dual is related to diagonal 
Yoneda Ext-algebra via
\[
\overline{A}^{!} \cong \bigoplus_i \Ext^i_{\overline{A}}(\Q,\Q)_i.
\]}
\end{enumerate}
Results of Peeva \cite{P} and Roos \cite{R} show that in general there 
does not exist a standard graded algebra $R$ such that 
$\prod_{k=1}^{\infty}(1-t^k)^{\phi_k}=HS(R,-t)$. 
For any quotient of a free Lie algebra, can we:
\begin{prob} Find spaces for which there is a simple generating function for
$\phi_k$.
\end{prob}
 \begin{prob} Relate $\hA$ to $\bigoplus\limits_{X \in L_2} \!\! {\mathfrak{h}}_{{\A}_X}$, as in \cite{PS}.
\end{prob}
As Shelton-Yuzvinsky proved in \cite{SY}, the natural class of arrangements for 
which an LCS formula holds are arrangements for which $A$ is
a Koszul algebra, which we tackle next.
\section{Koszul algebras}
Let $T(V)$ denote the tensor algebra on $V$.
\begin{defn}
A quadratic algebra is $T(V)/I$, with $I \subseteq V\otimes V$.
\end{defn}
A quadratic algebra $A$ has a quadratic dual $A^{\perp}=T(V^{*})/I^{\perp}$:
\[
\langle \alpha\otimes\beta \mid \alpha(a)\cdot \beta(b)=0 \mid \forall a\otimes b \in I \rangle = I^{\perp} \subseteq V^{*}\otimes  V^{*}
\]
\begin{defn} $A$ is Koszul if $Tor^A_i(\K,\K)_j = 0$, $i\ne j$.
\end{defn}
A quadratic algebra $A$ is Koszul iff the minimal free resolution of the residue
field over $A$ has matrices with only linear entries.  
\begin{exm}
The Hilbert series of 
$S = T(\K^n)/\langle x_i\otimes x_j - x_j\otimes x_i \rangle$ is $1/(1-t)^n$, and a computation
shows that the minimal free resolution of $\K$ over $S$ is the Koszul complex, so 
$\dim_{\K}Tor^S_i(\K,\K)_i = {n \choose i}$. Since
\[
I^{\perp} = \langle x_i\otimes x_j + x_j\otimes x_i \rangle, 
\]
we see that $S^{!}=E$. The Hilbert series of $E$ is $(1+t)^n = \sum_{i=0}^n {n \choose i}t^i$.
A computation shows that $\dim_{\K}Tor^E_i(\K,\K)_i = {n-1+i \choose i}$, which are
the coefficients in an expansion of $1/(1-t)^n$.
\end{exm}
Fr\"oberg \cite{Froberg} proved that if $I$ is a quadratic monomial ideal then $S/I$ is
Koszul. By uppersemicontinuity \cite{Herzog}, this means $S/I$ is Koszul if $I$ has a quadratic 
Grobner basis (QGB). See Example~\ref{KnoQGB} below for a Koszul algebra having
no QGB. Both $S$ and $E$ are Koszul, and the relation between their 
Hilbert series is explained by:
\begin{thm}
If $A$ is Koszul, so is $A^{!}$, and 
\[
HS(A,t)\cdot HS(A^{!}, -t) = 1
\]
\end{thm}
\begin{thm}[Bj\"{o}rner-Ziegler \cite{BZ}] The Orlik-Solomon algebra has a QGB iff $\A$ is supersolvable.
\end{thm}
\begin{exm}
A computation shows that the Orlik-Solomon algebra of $A_3$ has a quadratic Grobner basis,
so is Koszul. For the non-Fano arrangement, $\dim_{\K}Tor^A_3(\K,\K)_4=1$, so $A$ is not Koszul.
\end{exm}
\begin{exm}[Caviglia \cite{Cav}]\label{KnoQGB}
Map $R=\K[a_1,\ldots a_9] \stackrel{\phi}{\longrightarrow} \K[x,y,z]$ using all 
cubic monomials of $\K[x,y,z]$ except $xyz$, and let $I=\ker(\phi)$. Then $R/I$ 
is Koszul, but has no quadratic Grobner basis.
\end{exm}

\begin{prob} For Orlik-Solomon algebras, does Koszul imply supersolvable? In the
case of graphic arrangements, it does \cite{SS}.
\end{prob}
\begin{prob} Find a combinatorial description of $Tor^A_i(\K,\K)_j$.
\end{prob}
\section{Resonance varieties}\label{resonanceSection}
Let $A$ be the Orlik-Solomon algebra of $\XA$, with $|\A|=n$. 
For each $a=\sum a_ie_i \in A_1$, we consider the Aomoto complex $(A,a)$,
whose $i^{\text{th}}$ term is $A_i$, and differential is $ \wedge a$:
\[
(A,a)\colon 
\xymatrix{
0 \ar[r] &A_0 \ar[r]^{a} & A_1
\ar[r]^{a}  & A_2 \ar[r]^{a}& \cdots \ar[r]^{a}
& A_{\ell}\ar[r] & 0}.
\]
This complex arose in Aomoto's work \cite{Ao} on 
hypergeometric functions, as well as in the study of 
cohomology with local system coefficients \cite{ESV}, \cite{STV}.
In \cite{Yuz}, Yuzvinsky showed that for a generic $a$, the
Aomoto complex is exact; the resonance varieties of 
$\A$ are the loci of points
$a=\sum_{i=1}^na_ie_i \leftrightarrow (a_1:\dots :a_n) \in \P^{n-1}$ for 
which $(A,a)$ fails to be exact, that is:
\vskip -.2in
\begin{defn} For each $k\ge 1$,
\[
R^{k}(\A)=\{a \in \P^{n-1} \mid  H^k(A,a)\ne 0\}.
\]
\end{defn}
In \cite{Fa}, Falk gave necessary and sufficient
conditions for $a \in R^1(\A)$.
 \begin{defn} A  partition $\Pi$ of $\A$ is neighborly if 
for all $Y \in L_2(\A)$ and $\pi$ a block of $\Pi$,
\[
\mu(Y) \le  |Y\cap \pi| \Longrightarrow Y\subseteq \pi.
\]
\end{defn}
Falk proved that components of $R^1(\A)$ arise from 
neighborly partitions; he conjectured that $R^1(\A)$ is a union 
of linear components. This was established (essentially simultaneously)
by Cohen--Suciu \cite{CScv} and Libgober-Yuzvinsky \cite{LY}. 
Libgober and Yuzvinsky also showed that $R^{1}(\A)$ is a disjoint
union of positive dimensional subspaces in $\P(E_1)$, and Cohen-Orlik 
\cite{CO} show that $R^{k \ge 2}(\A)$ is also a subspace arrangement. 

On the other hand, as shown by Falk in \cite{Fa2}, in positive characteristic,
the components of $R^{1}(\A)$ can meet, and need not be linear.
The approach of Libgober--Yuzvinsky involves connecting $R^1(\A)$ to 
pencils/nets/webs and there is much recent work in the area, e.g.
\cite{FalkY}, \cite{PY} \cite{yu3}. Of special interest is the following 
conjecture relating $R^1(\A)$ and the LCS ranks $\phi_k$:
\begin{conj}[Suciu \cite{Su}]
If $\phi_4 = \theta_4$, then 
\[
\prod\limits_{k \ge 1} (1-t^k)^{\phi_k} = \prod\limits_{L_i \in
    R^1(\A)}(1-(dim(L_i)t),
\]
where $\theta_4$ is the fourth Chen rank (Definition~\ref{Cranks}).
\end{conj}
  
\begin{exm}\label{exm:toyRes}
Let $\A = V(xy(x-y)z) \subseteq \P^2$, and 
$E = \Lambda(\K^4)$, with generators $e_1,\ldots, e_4$, so that 
\[
A = E/\langle \partial(e_1e_2e_3), \partial(e_1e_2e_3e_4) \rangle.
\]
Since $\partial(e_1e_2e_3e_4)  =  e_1 \wedge \partial(e_1e_2e_3) - e_4 \partial(e_1e_2e_3)$,
the second relation is unnecessary. To compute $R^1(\A)$, we need only the first two
differentials in the Aomoto complex. Using 
$e_{13},e_{14},e_{23},e_{24},e_{34}$ as a basis for $A_2$, we find that 
$e_1 \mapsto e_1 \wedge (\sum\limits_{i=1}^4 a_ie_i) = a_2 e_{12} + a_3 e_{13} + a_4 e_{14}$.
Since 
\[
\partial(e_1e_2e_3)  =  e_1 \wedge e_2 -e_1 \wedge e_3+ e_2 \wedge e_3,
\]
in $A$, $e_{12} = e_{13}-e_{23}$, so that 
\[a_2 e_{12}=a_2(e_{13}-e_{23}).
\]
This means $e_1 \mapsto (a_2 + a_3)e_{13} + a_4 e_{14}-a_2e_{23}$;
similar computations for the other $e_i$ show that the Aomoto 
complex is
\[
0 \longrightarrow \K^1 \xrightarrow{\left[ \!
\begin{array}{c}
a_1\\
a_2\\
a_3\\
a_4
\end{array}\! \right]} \K^4
\xrightarrow{\left[ \!\begin{array}{cccc}
a_2+a_3 & -a_1     & -a_1 & 0\\
a_4     & 0        & 0 & -a_1\\
-a_2    &  a_1+a_3 & -a_2 & 0 \\
0       & a_4      & 0 & -a_2 \\
0       & 0        & a_4 & -a_3
\end{array}\! \right]}
\K^5
\]
The rank of the first map is always one, $R^1(A)\subseteq \P^3$ is the locus where the second 
matrix has kernel of dimension at least two, so the $3 \times 3$ minors must vanish.
A computation shows this locus is $\langle a_4, a_1+a_2+a_3\rangle$. 
\end{exm}
Letting $a=\sum_{i=1}^na_ie_i$, observe that $a \in R^1(\A)$ iff there exists a
$b \in E_1$ so that $a\wedge b \in I_2$, so that $R^1(\A)$ is the locus of decomposable 
2-tensors in $I_2$. Since $I_2$ is determined by the intersection lattice $L(\A)$ 
in rank $\le 2$, to study $R^1(\A)$, it can be assumed that $\A \subseteq \P^2$.

While the first resonance variety is conjecturally connected (under certain conditions) 
to the LCS ranks, $R^1(A)$ is {\em always} connected to the Chen ranks introduced by 
Chen in \cite{C}.
\begin{defn}\label{Cranks}
The Chen ranks of $G$ are the LCS ranks of the maximal 
metabelian quotient of $G$:
\[
\theta_k(G):=\phi_k(G/G''),
\] 
where $G'=[G,G]$. 
\end{defn}
\pagebreak
\begin{conj}[Suciu \cite{Su}]
\label{conj:chen}
Let $G=G(\A)$ be an arrangement group, and let 
$h_r$ be the number of components
of $R^{1}(\A)$ of dimension $r$.  Then, for $k \gg 0$:
\[
\theta_k(G)= (k-1) \sum_{r\ge 1} h_r  \binom{r+k-1}{k}.  
\]
\end{conj}
For the previous example, $R^1(\A) = V(a_4, a_1+a_2+a_3) \simeq \P^1$, so
\[
\theta_k(G)= (k-1).  
\]
To discuss the Chen ranks, we need some background. 
The Alexander invariant $G'/G''$ is a module over $\Z[G/G']$. For 
arrangements, $\Z[G/G']=$ Laurent polynomials in $n$-variables.
In \cite{Massey}, Massey showed that 
\[
\sum_{k\ge 0} \theta_{k+2}\, t^k = HS(\mbox{gr }G'/G'' \otimes \Q,t). 
\]
It turns out to be easier to work with the linearized Alexander invariant $B$ 
introduced by Cohen-Suciu in \cite{CS2}
\[
(A_2 \oplus E_3)\otimes S \stackrel{\Delta}{\longrightarrow} 
E_2\otimes S \longrightarrow B \longrightarrow  0,
\]
where $\Delta$ is built from the Koszul differential and 
$(E_2\rightarrow A_2)^t$.
\begin{thm}[Cohen-Suciu \cite{CS2}] 
\[
V(\mbox{ann }B)=R^1(\A)
\]
\end{thm}

\begin{thm}[Papadima-Suciu \cite{PSChen}] For $k\ge 2$,
\[
\sum_{k\ge 2} \theta_{k}\, t^k = HS(B,t).
\]
\end{thm}
This shows that the Chen ranks are combinatorially determined,
and depend only on $L(\A)$ in rank $\le 2$.
\begin{exm} 
For the $A_3$ arrangement depicted in Example~\ref{ShowMeThm},
write $e_0=L_{12}, e_1=L_{13},e_2=L_{23},e_3=L_{24}, e_4=L_{14},e_5=L_{34}$.
With this labelling 
\[
I_2 = \langle \partial(e_1e_4e_5), \partial(e_0e_1e_2),\partial(e_2e_3e_5), \partial(e_0e_3e_4)\rangle,
\]
from which a presentation for $B$ can be computed: 
\[
S^{14} \rightarrow S^4 \rightarrow B \rightarrow 0.
\]
A computation shows that $R^1(A_3)$ is
\begin{center}
$\begin{array}{c}
V(x_1 + x_4 + x_5,x_0,x_2,x_3)\sqcup\\ 
V(x_2 + x_3 + x_5,x_0,x_1,x_4)\sqcup\\
V(x_0 + x_3 + x_4,x_1,x_2,x_4)\sqcup\\
V(x_0 + x_1 + x_2,x_3,x_4,x_5)\sqcup\\
V(x_0+x_1+x_2,x_0- x_5,x_1-x_3,x_2-x_4),
\end{array}$
\end{center}
and that the Hilbert Series of $B$ is:
\[
(4t^2+2t^3-t^4)/(1-t)^2 = 4t^2+10t^3+15t^4+20t^5+\cdots 
\]
On the other hand, the graded betti numbers 
$Tor^E_i(A_3,\K)_j$ are
$$
\vbox{\offinterlineskip 
\halign{\strut\hfil# \ \vrule\quad&# \ &# \ &# \ &# \ &# \ &# \ &# \ &# \ &# \ &# \ &#\ &# \ &#\  \cr
total&1 &4 &10&21&45&91   \cr \noalign {\hrule} 
0    &1 &--&--&--&-- &--   \cr 
1    &--&4 &10&15&20 &25   \cr 
2    &--&--&--&6& 25& 66  \cr 
}}
$$
So the Hilbert series for $B$ encodes the ranks of $Tor^E_i(A_3,\K)_{i+1}$.
This suggests a connection between $R^1(A)$ and $Tor^E_i(A_3,\K)_{i+1}$,
which we tackle in the next section.
\end{exm}
Besides the connection to resonance varieties, there is a second 
reason to study $Tor^E_i(A,\K)$: the numbers 
$b_{ij} = \dim_{\K}Tor^A_i(\K,\K)_j$ studied in \S 7 grow very fast,
while the numbers $b_{ij}' = \dim_{\K}Tor^E_i(A,\K)_j$ grow at a 
much slower rate. 
\begin{exm}
For the non-Fano arrangement of Example~\ref{exm:nonfano}\newline
$$
\vbox{\offinterlineskip 
\halign{\strut\hfil# \ \vrule\quad&# \ &# \ &# \ &# \ &# \ &# \ &# \ &# \ &# \ &# \ &#\ &# \ &#\  \cr
total&1 &7 &23&63&165&387   \cr \noalign {\hrule} 
0    &1 &--&--&--&-- &--   \cr 
1    &--&6 &17&27&36 &45   \cr 
2    &--&1 &6 &36&129&342  \cr 
3    &-- &--&--&--&--&--  \cr 
\noalign{\bigskip} \noalign{\smallskip} }}
\vbox{\offinterlineskip 
\halign{\strut\hfil# \ \vrule\quad&# \ &# \ &# \ &# \ &# \ &# \ &# \ &# \ &# \ &# \ &#\ &# \ &# \ &#\  \cr
total&1&7 &35& 156& 664& 2773 &*\cr \noalign {\hrule} 
0&1 &7  &34   & 143&560 & 2108&* \cr 
1&--&-- &1    & 13&103 &646 &* \cr 
2&--&-- &--   & --&1  & 19 &* \cr 
3&--&-- &--   &--&--  &-- &1  \cr 
\noalign{\bigskip} \noalign{\smallskip} }}
$$
\vskip -.2in
\hskip 1in $b_{ij}'$  \hskip 2in $b_{ij}$.
\end{exm}
The spaces $Tor^E_i(A,\K)$ and $Tor^A_i(\K,\K)$ are related via
the change of rings spectral sequence:
\[
\Tor_i^A\left(\Tor_j^E(A,\K),\K\right)
\Longrightarrow \Tor_{i+j}^E(\K,\K).
\]
For arrangements, details of this relationship are investigated in
 \cite{SS}.
\begin{prob} Find a combinatorial description of $Tor^E_i(A,\K)_j$.
\end{prob}  
\begin{prob} If $A$ is Koszul, does this provide data on $Tor^E_i(A,\K)_j$?
\end{prob}

\section{Linear syzygies}
It is fairly easy to see that there is a connection between $R^1(A)$ and linear 
syzygies, that is, to the module $Tor^E_2(A,\K)_3$. Since
\[
a\wedge b \in I_2 \longrightarrow  a\wedge b = \sum c_if_i, \mbox{ }c_i \in \K, f_i \in I_2,
\]
the relations $a\wedge a\wedge b = 0 = b\wedge a\wedge b$ yield linear syzygies on $I_2$:
\[
\sum ac_if_i = 0 =\sum bc_if_i.
\]
{\bf Example~\ref{exm:toyRes}}, continued. Since $\partial(e_1e_2e_3) = (e_1-e_2) \wedge (e_2-e_3)$,
both $(e_1-e_2)$ and $(e_2-e_3)$ are in $R^1(A)$, as is the line connecting them:
\[
s(e_1-e_2)+ t(e_2-e_3) \subseteq R^1(\A) \subseteq \P(E_1)
\]
Parametrically, this may be written
\[
(s : t-s : -t:0) = V(a_4, a_1+a_2+a_3),
\]
so $s(e_1-e_2)+ t(e_2-e_3) \wedge \partial(e_1e_2e_3) = 0$ 
gives a family of linear syzygies on $I_2$, 
parameterized by $\P^1$. $\Diamond$\newline

To make the connection between linear syzygies and the module $B$ precise,
we need the following result:
\begin{thm}[Eisenbud-Popescu-Yuzvinsky \cite{EPY}] For
an arrangement $\A$, the Aomoto complex is exact,
as a sequence of $S$-modules:
\[
\xymatrixcolsep{21pt}
\xymatrix{
0\ar[r] & A_0  \otimes S \ar[r]^{\cdot a} &
A_1\otimes S \ar[r]^(.62){\cdot a} & \cdots \ar[r]^(.37){\cdot a}  & 
A_{\ell} \otimes S \ar[r] & F(A) \ar[r] & 0
}.
\]
\end{thm}

\begin{thm}[Schenck-Suciu \cite{SS2}]
The linearized Alexander invariant $B$ is determined by $F(A)$:
\[
B \cong\Ext_S^{\ell-1}(F(A), S).
\]
Furthermore,  for $k\ge 2$, $\dim_{\K}B_k = \dim_{\K}\Tor^E_{k-1}(A,\K)_k.$
\end{thm}
Using this, it is possible to prove one direction of Conjecture~\ref{conj:chen}
\begin{thm}[Schenck-Suciu \cite{SS2}]
For $k \gg 0$,
\[
\theta_k(G) \ge (k-1) \sum_{L_i \in R^1(\A)}  \binom{\dim L_i+k-1}{k}.
\]
\end{thm}
\begin{prob} Prove the remaining direction of Conjecture~\ref{conj:chen}.
\end{prob}
What makes all this work is the Bernstein-Gelfand-Gelfand correspondence, which is
our final topic.
\pagebreak

\section{Bernstein-Gelfand-Gelfand correspondence}
Let $S=Sym(V^*)$ and $E = \bigwedge(V)$. The 
Bernstein-Gelfand-Gelfand correspondence is an 
isomorphism between derived categories of bounded complexes 
of coherent sheaves on $\P(V^*)$ and bounded complexes of 
finitely generated, graded $E$--modules. Although this
sounds exotic, from this it is possible to extract
functors\newline

\noindent ${\bf R}$: finitely generated, graded $S$-modules $\longrightarrow$ linear free $E$-complexes.
\newline
\noindent${\bf L}$: finitely generated, graded $E$-modules $\longrightarrow$ linear free $S$-complexes.\newline

The point is that problems can be translated to a (possibly) simpler setting. 
For example, BGG yields a very fast way to compute sheaf cohomology, using Tate 
resolutions.

\begin{defn}
Let $P$ be a finitely generated, graded $E$-module. Then
${\bf L}(P)$ is the complex
\[
\xymatrixcolsep{21pt}
\xymatrix{
 \cdots\ar[r] &S \otimes P_{i+1}   \ar[r]^(.56){\cdot a} &
S \otimes P_{i}   \ar[r]^(.44){\cdot a} &S \otimes P_{i-1}  \ar[r]^(.66){\cdot a} &\cdots
},
\]
where $a=\sum\limits_{i=1}^n x_i\otimes e_i$, so that 
$1 \otimes p \mapsto \sum x_i\otimes e_i \wedge p$
\end{defn}
Note that elements of $V^*$ have degree $=1$, and elements of $V$ have degree $=-1$.
\begin{exm}
$P = E = \bigwedge \K^3$. Then we have
\[
0 \longrightarrow S \otimes E_0 \longrightarrow S \otimes E_1
\longrightarrow  S \otimes E_2 \longrightarrow  S \otimes E_3 \longrightarrow 0.
\]
Clearly $1 \mapsto \sum_{1}^3 x_i\otimes e_i$. For $d_1$ 
\begin{center}
$\begin{array}{c}
\;\;\;e_1 \mapsto -x_2e_{12}-x_3e_{13}\\
e_2 \mapsto x_1e_{12}-x_3e_{23}\\
e_3 \mapsto x_1e_{13}+x_2e_{23}
\end{array}$
\end{center}
\[
d_2:\mbox{ }e_{12} \mapsto x_3 e_{123}, \mbox{ }e_{13}\mapsto -x_2 e_{123}
\mbox{ }e_{23}\mapsto x_1 e_{123}
\]
Thus,  ${\bf L}(E)$ is
\begin{small}
\[
S^1 \xrightarrow{\left[ \!
\begin{array}{c}
x_1\\
x_2\\
x_3
\end{array}\! \right]} S^3
\xrightarrow{\left[ \!\begin{array}{ccc}
-x_2 & x_1   & 0 \\
-x_3 & 0     & x_1 \\
 0   & -x_3  & x_2
\end{array}\! \right]}
S^3
\xrightarrow{\left[ \!\begin{array}{ccc}
x_3  & -x_2 & x_1
\end{array}\! \right]}
S^1
\]
\end{small}
This is simply the Koszul complex.
\end{exm}
If $M$ is a finitely generated, graded $S$-module, then
${\bf R}(M)$ is the complex
\[
\xymatrixcolsep{21pt}
\xymatrix{
 \cdots\ar[r] &\hat{E} \otimes M_{i-1}   \ar[r]^(.56){\cdot a} &
\hat{E} \otimes M_{i}   \ar[r]^(.44){\cdot a} &\hat{E} \otimes M_{i+1} \ar[r]^(.66){\cdot a} &\cdots
},
\]
where $a=\sum\limits_{i=1}^n e_i\otimes x_i$, so 
$1 \otimes m \mapsto \sum e_i\otimes x_i \cdot m$,
and $\hat{E}$ is $\K$-dual to $E$:
\[
\hat{E} \simeq E(n)=Hom_{\K}(E,\K).
\]
Just as ${\bf L}(P) = S\otimes_\K P$, 
${\bf R}(M) = Hom_{\K}(E,M)$.
\vskip .05in
{\bf Example~\ref{exm:first}}, continued.
If $M = \K[x_0,x_1]/\langle x_0x_1, x_0^2 \rangle$, then 
\begin{center}
$\begin{array}{c}
1 \mapsto e_0\otimes x_0 + e_1\otimes x_1 \\
\;\;x_0 \mapsto e_0\otimes x_0^2 +e_1\otimes x_0x_1 \\
\;\;x_1 \mapsto  e_0\otimes x_0x_1 +e_1\otimes x_1^2\\
\;\;\;\;\;\;x_1^n \mapsto  e_0\otimes x_0x_1^n +e_1\otimes x_1^{n+1}
\end{array}$
\end{center}
Thus,  ${\bf R}(M)$ is
\begin{small}
\[
E(2)^1 \xrightarrow{\left[ \!
\begin{array}{c}
e_0\\
e_1
\end{array}\! \right]} E(3)^2
\xrightarrow{\left[ \!\begin{array}{cc}
 0  & e_1
\end{array}\! \right]}
E(4)^1
\xrightarrow{\left[ \!\begin{array}{c}
e_1
\end{array}\! \right]}
E(5)^1
\xrightarrow{\left[ \!\begin{array}{c}
e_1
\end{array}\! \right]}
 \cdots
\]
\end{small}
This complex is exact, except at the second step. The kernel of 
\[
\left[ \!\begin{array}{cc}
 0  & e_1
\end{array}\! \right]
\]
is generated by $\alpha = [1,0]$ and $\beta = [0,e_1]$, with
relations $im(d_1)=\beta+e_0 \alpha =0, e_1 \beta = 0$, so that 
\[
H^1({\bf R}(M))  \simeq E(3)/e_0 \wedge e_1
\]
The betti table for $M$ is:
$$
\vbox{\offinterlineskip 
\halign{\strut\hfil# \ \vrule\quad&# \ &# \ &# \ &# \ &# \ &# \
&# \ &# \ &# \ &# \ &# \ &# \ &# \
\cr
total&1&2&1\cr
\noalign {\hrule}
0&1 &--&--&\cr
1&--&2 &1 &\cr
\noalign{\bigskip}
\noalign{\smallskip}
}}
$$
\vskip -.2in
Note that in this example, $M$ is $1$-regular.$\Diamond$ 

\begin{thm}[Eisenbud-Fl{\o}ystad-Schreyer \cite{EFS}]
\[
H^j({\bf R}(M))_{i+j} = Tor_i^S(M,\K)_{i+j}.\]
\end{thm}
\begin{corollary}\label{EFSCM}
The Castelnuovo-Mumford regularity of $M$ is $\le d$ iff
$H^i({\bf R}(M))=0$ for all $i > d$.
\end{corollary}
\pagebreak 
What can be said about higher resonance varieties? In \cite{CO}, Cohen--Orlik prove 
that for $k \ge 2$,
\[
R^k(\A) = \bigcup L_i \mbox{ linear.}
\]
As observed by Suciu, in general the union need not be disjoint. 
\begin{thm}[Eisenbud-Popescu-Yuzvinsky \cite{EPY}] $R^k(\A) \subseteq R^{k+1}(\A).$
\end{thm}
The key point is that
\[
H^k(A,a) \ne 0 \mbox{ iff } Tor^S_{\ell-k}(F(A), S/I(p)) \ne 0.
\]
The result follows from interpreting this in terms of Koszul cohomology.
\begin{thm}[Denham-Schenck \cite{denS}]
Higher resonance may be interpreted via $\Ext$:
\[
R^k(\A) = \bigcup_{k' \le k} V(\ann \Ext^{\ell-k'}(F(A),S)).
\]
Furthermore, the differentials in free resolution of $A$ over $E$ can be analyzed 
using BGG and the Grothendieck spectral sequence.
\end{thm}
For any coherent sheaf ${\mathcal F}$ on $\P^d$, there is a finitely generated,
graded {\em saturated} $S$-module $M$ whose sheafification is ${\mathcal F}$.
If ${\mathcal F}$ has Castelnuovo-Mumford regularity $r$, then 
the Tate resolution of ${\mathcal F}$ is obtained by 
splicing the complex ${\bf R}(M_{\ge r})$:
\[
\xymatrixcolsep{21pt}
\xymatrix{
0\ar[r] &\hat{E} \otimes M_{r}   \ar[r]^{d^r} &
\hat{E} \otimes M_{r+1}   \ar[r] &\hat{E} \otimes M_{r+2} \ar[r] &\cdots
},
\]
with a free resolution $P_{\bullet}$ for the kernel of $d^r$:
\[
\xymatrixrowsep{10pt}
\xymatrixcolsep{6pt}
\xymatrix{
\cdots \ar[rr] && P_1 \ar[rr] && P_0 \ar[rr] \ar[dr]  &                   & \hat{E} \otimes M_{r}  \ar[rr]&&  \hat{E} \otimes M_{r+1} \ar[rr]   && \cdots \\
               &&             &&                      & \ker(d^r) \ar[ur] \ar[dr] &            &&               &&        \\
               &&             && 0 \ar[ur]            &                   &  0         &&               &&      
}
\]

By Corollary~\ref{EFSCM}, ${\bf R}(M_{\ge r})$ is exact except at the first step, so this yields an exact complex of free $E$-modules. 
\begin{exm}
Since $M=S$ has regularity zero, we obtain Cartan resolutions in both directions, and the
splice map $ E \rightarrow \widehat{E}=E(d+1)$ is multiplication by $e_0 \wedge e_1 \wedge \cdots \wedge e_d = \ker \left[ \!
\begin{array}{ccc}
e_0, & \cdots &, e_d \end{array}\! \right]^t$.
\end{exm}
\pagebreak
\begin{thm}[Eisenbud-Fl{\o}ystad-Schreyer \cite{EFS}]\label{thmEFS}
The $i^{th}$ free module $T^i$ in a Tate resolution for ${\mathcal F}$
satisfies
\begin{center}
$T^i = \bigoplus\limits_{j}\widehat{E} \otimes H^j({\mathcal F}(i-j))$.\end{center}
\end{thm}

{\bf Example~\ref{exm:twistedcubic}}, continued.
The betti table for the twisted cubic shows that $S/I$ has regularity one, 
which provides us the information needed to compute the Tate resolution. 
Plugging the resulting numbers into Theorem~\ref{thmEFS} shows that

\begin{center}
\begin{supertabular}{|c|c|c|c|c|c|c|} 
\hline $i$                    & $-3$ & $-2$ & $-1$ &  $0$ & $1$ & $2$ \\
\hline $h^1({\mathcal F}(i))$ & $8$ & $5$  & $2$ & $0$ & $0$ & $0$\\
\hline $h^0({\mathcal F}(i))$ & $0$ & $0$  & $0$ & $1$ & $4$ & $7$  \\
\hline
\end{supertabular}
\end{center}
\vskip .05in
Does this make sense? Since ${\mathcal F} = {\mathcal O}_X = {\mathcal O}_{\P^1}(3)$,
\[
h^1({\mathcal F}(i)) = h^1({\mathcal O}_{\P^1}(3i))=h^0({\mathcal O}_{\P^1}(-3i-2))
\]
and 
\[
h^0({\mathcal F}(i)) = h^0({\mathcal O}_{\P^1}(3i))=3i+1, \mbox{ }i\ge 0,
\]
which agrees with our earlier computation. $\Diamond$

\vskip .05in
\begin{prob}
Investigate the Tate resolution for $D(\A)$ and $C(\A)$. 
\end{prob}
\noindent{\bf Conclusion} In this note we have surveyed a number
of open problems in arrangements. The beauty of the area is that these
problems are all interconnected. Perhaps the most central objects are
the resonance varieties, which are related to both the LCS ranks studied
in \S 7 and \S 8 using Koszul and Lie algebras, and to the Chen ranks. 
The results 
of \S 5 tie resonance to the Orlik-Terao algebra, and \cite{syzRes}
notes that $J_{\A} \subseteq H^0(D_{\A})$, so the Orlik-Terao 
algebra is also linked 
to $D(\A)$ and freeness. But freeness ties in to multiarrangements, and 
can be generalized to hypersurface arrangements, the topics of \S 3 and \S 4.
To complete the circle, recent work of 
Cohen-Denham-Falk-Varchenko \cite{CDFV} relates
freeness to $R^1(\A)$. In short, everything is connected!
\vskip .05in
\noindent{\bf Acknowledgements} Many thanks are due to the organizers of the Mathematical
Society of Japan Seasonal Institute: Takuro Abe, Hiroaki Terao, Masahiko Yoshinaga and 
Sergey Yuzvinsky organized a wonderful meeting. I also am grateful to my research 
collaborators: G.~Denham, M.~Mustata, A.~Suciu, H.~Terao, S.~Tohaneanu, and M.~Yoshinaga.
As noted in the introduction, the calculations carried out in this survey
may be performed using the hyperplane arrangements package of 
Denham and Smith~\cite{DSmith} in Macaulay2~\cite{GS}.
Code for the individual examples is available at 
{\tt http://www.math.uiuc.edu/~schenck/}.



\begin{thebibliography}{99}

\bibitem{atw} T.~Abe, H.~Terao, M.~Wakefield,
{\em The Euler multiplicity and addition-deletion theorems for multiarrangements}, 
J. Lond. Math. Soc. \textbf{77} (2008), 335--348.

\bibitem{atw2} T.~Abe, H.~Terao, M.~Wakefield,
{\em The characteristic polynomial of a multiarrangement}, 
Adv. in Math,  \textbf{215}  (2007), 825--838. 

\bibitem{aty}T.~Abe, H.~Terao, M.~Yoshinaga,
{\em Totally free arrangements of hyperplanes},
Proc. Amer. Math. Soc.  \textbf{137}  (2009),1405--1410. 

\bibitem{Ao} K.~Aomoto,
{\em Un th\'eor\`eme du type de Matsushima-Murakami concernant l'int\'egrale des fonctions multiformes},
J. Mat. Pur. App. \textbf{52}  (1973), 1--11.

\bibitem{Ar} W.~Arvola, 
{\em The fundamental group of the complement of an arrangement of complex hyperplanes}, 
    Topology  \textbf{31}  (1992), 757--765.

\bibitem{BZ} A.~Bj\"{o}rner, G.~Ziegler,
{\em Broken circuit complexes: factorizations and 
generalizations},  J. Combin. Theory Ser. B \textbf{51} 
(1991), 96--126.

\bibitem{bt} K. ~Brandt, H. ~Terao,
        {\em Free Arrangements and Relation Spaces},
        Discrete Comput. Geom., \textbf{12} (1994), 49-63.

\bibitem{Br} E.~Brieskorn,
              {\em Sur les groupes de tresses}, 
              S\'{e}minaire Bourbaki, 1971/72, 
              LNM \textbf{317}, Springer-Verlag, (1973), 21--44. 

\bibitem{Cav} G.~Caviglia,
            {\em The pinched veronese is Koszul},
            J. Algebraic Combin.,  \textbf{30}, (2009), 539--548.

\bibitem{C} K.~Chen,
            {\em Iterated integrals of differential forms and loop space cohomology},
              Ann. of Math. \textbf{97} (1973), 217-246.

\bibitem{CDFV}  D.~Cohen, G.~Denham, M.~Falk, A.~Varchenko,
{\em Critical points and resonance of hyperplane arrangements}, 
Canadian J. Math, \textbf{63}, (2011), 1038-1057.

\bibitem{CO}  D.~Cohen, P.~Orlik,
{\em Arrangements and local systems}, Math. Res. Lett.
\textbf{7} (2000), 299--316.

\bibitem{CS} D.~Cohen, A.~Suciu,
 {\em The braid monodromy of plane algebraic curves and hyperplane arrangements},
Comment. Math. Helv.  \textbf{72} (1997), 285--315.       

\bibitem{CScv} D.~Cohen, A.~Suciu,
{\em Characteristic varieties of arrangements},
Math. Proc. Cambridge Phil. Soc. \textbf{127} (1999), 33--53.

\bibitem{CS2} D.~Cohen, A.~Suciu,
             {\em Alexander invariants of complex hyperplane arrangements},
             Trans. Amer. Math. Soc.  \textbf{351} (1999), 4043--4067.

\bibitem{csTapp} D.~Cohen, A.~Suciu, eds.
             {\em Arrangements in Boston: A conference on hyperplane arrangements}
              Topology Appl.  \textbf{118}  (2002)

\bibitem{DE} W.~Decker, D.~Eisenbud,
{\em Sheaf algorithms using the exterior algebra}, in:
Computations in Algebraic Geometry using Macaulay~2, 
Springer-Verlag, Berlin-Heidelberg-New York, 2002. 

\bibitem{DP} C.~De Concini, C.~Procesi,
{\em Wonderful models of subspace arrangements}
 Selecta Math. \textbf{1}  (1995),  no. 3, 459--494.

\bibitem{DP2} C.~De Concini, C.~Procesi, 
        {\em Hyperplane arrangements and holonomy equations},
        Sel. Math. {\bf 1} (1995), 495--535.

\bibitem{den} G.~Denham, 
{\em Homological aspects of hyperplane arrangements}, in:
Arrangements, Local Systems and Singularities, Birkhauser 2010.

\bibitem{denS} G.~Denham,  H.~Schenck,
{\em The double $\Ext$ spectral sequence and rank varieties},
preprint, 2010.

\bibitem{DSmith} G.~Denham,  G.G.~Smith,
{\em Hyperplane Arrangements: a Macaulay2 package}
{\tt http://www.math.uiuc.edu/Macaulay2/}

\bibitem{DS} G.~Denham, A.~Suciu,
             {\em On the homotopy Lie algebra of an arrangement},
             Michigan Math. J.  \textbf{54}  (2006),  319--340.

\bibitem{DY} G.~Denham,  S.~Yuzvinsky,
{\em Annihilators of Orlik-Solomon relations},
 Adv. in Appl. Math. \textbf{28}  (2002),  231--249.

\bibitem{dimcaYuz} A.~Dimca, S.~Yuzvinsky, 
{\em Lectures on Orlik-Solomon algebras}, in:
Arrangements, Local Systems and Singularities, Birkhauser 2010.

\bibitem{EisenbudFactor} D. ~Eisenbud,
        {\em Linear sections of determinantal varieties},
        Amer. J. Math.  {\bf 110}  (1988),  no. 3, 541--575. 

\bibitem{E} D. ~Eisenbud,
{\em The geometry of syzygies}, 
Springer-Verlag, Berlin-Heidelberg-New York, 2004. 

\bibitem{EFS} D.~Eisenbud,  G.~Fl{\o}ystad, F.-O.~Schreyer,
{\em Sheaf cohomology and free resolutions over exterior algebras},
Trans. Amer. Math. Soc. \textbf{355}  (2003),  4397--4426.  

\bibitem{EPY} D.~Eisenbud,  S.~Popescu, S.~Yuzvinsky,
{\em Hyperplane arrangement cohomology and monomials in the exterior algebra}, 
Trans. Amer. Math. Soc. \textbf{355}  (2003), 4365--4383.  

\bibitem{ESTUYV} F.~El Zein, A.~Suciu, M.~Tosun, M~Uludag, S.~Yuzvinsky (eds),
{\em Arrangements, Local Systems and Singularities}
Birkhauser, Progress in Mathematics, Vol. 283, 2010

\bibitem{ESV} H.~Esnault, V. ~Schechtman, E. ~Viehweg,
{\em Cohomology of local systems on the complement of hyperplanes},
Invent. Math. \textbf{109} (1992), 557-561.

\bibitem{f} M. ~Falk,
            {\em The minimal model of the complement of an arrangement of hyperplanes},
            Trans. Amer. Math. Soc. \textbf{309} (1988), 543--556.

\bibitem{Fa89} M. ~Falk,
            {\em The cohomology and fundamental group of a hyperplane complement},
            Contemporary Math., vol.~90, Amer. Math. Soc, Providence, RI, 1989, pp. 55--72.

\bibitem{Fa} M.~Falk,
{\em Arrangements and cohomology},
Ann. Combin. \textbf{1} (1997), 135--157.

\bibitem{Fa2} M.~Falk,
{\em Resonance varieties over fields of positive characteristic},
IMRN  \textbf{3} (2007) 

\bibitem{FR85} M. Falk, R. Randell,
            {\em The lower central series of a fiber-type arrangement},
            Invent. Math. \textbf{82} (1985), 77--88.

\bibitem{fr} M.~Falk, R.~Randell,
            {\em On the homotopy theory of arrangements II},
            Adv. Stud. Pure. Math. \textbf{27} (2000), 93-125.

\bibitem{ft} M.~Falk, H.~Terao eds.,
            {\em Arrangements---Tokyo 1998. Proceedings of the Workshop on Mathematics Related with Arrangements of Hyperplanes, celebrating the 60th birthday of Peter Orlik}
             Mathematical Society of Japan, Tokyo, 2000.

\bibitem{FalkY} M.~Falk, S.~Yuzvinsky, 
{\em Multinets, resonance varieties, and pencils of plane curves},
Compos. Math.  \textbf{143} (2007), 1069--1088 

\bibitem{F} E.~Feichtner,
{\em De Concini-Procesi wonderful arrangement models: a discrete geometer's point of view.}
Combinatorial and computational geometry,  333--360, Math. Sci. Res. Inst. Publ. 52, Cambridge Univ. Press, Cambridge, 2005. 

\bibitem{FK} E.~Feichtner, D.~Kozlov,
{\em Incidence combinatorics of resolutions},
Selecta Math. \textbf{10}  (2004),  no. 1, 37--60. 

\bibitem{FY} E.~Feichtner, S.~Yuzvinsky,
{\em Chow rings of toric varieties defined by atomic lattices.} 
Invent. Math. \textbf{155} (2004), no. 3, 515--536.

\bibitem{Froberg} R.~Fr\"{o}berg, 
{\em Determination of a class of Poincar\'e series}, 
Math. Scand. \textbf{37} (1975), 29--39.

\bibitem{FM} W.~Fulton, R.~MacPherson,
{\em A compactification of configuration spaces.}
Ann. of Math. {\bf 139} (1994), 183--225.

\bibitem{GS} D.~Grayson,  M.~Stillman,
{\em Macaulay2: a software package for commutative algebra}
{\tt http://www.math.uiuc.edu/Macaulay2/}

\bibitem{HL}  H.~Hamm, D.~T.~L\^e,
{\em Un th\'eor\`eme de Zariski du type de Lefschetz},
Ann. Sci. \'{E}cole Norm. Sup. \textbf{6} (1973), 317--366.

\bibitem{Herzog} J.~Herzog,
{\em Finite free resolutions},
 Computational commutative and non-commutative algebraic geometry,  118--144, NATO Sci. Ser. III 
Comput. Syst. Sci., 196, IOS, Amsterdam, 2005.

\bibitem{JP} M.~Jambu, S.~Papadima,
            {\em A generalization of fiber-type arrangements and a new deformation method},
            Topology \textbf{37} (1998), 1135--1164.

\bibitem{keel} S.~Keel,  
{\em Intersection theory of moduli space of stable n-pointed curves of genus zero},
Trans. Amer. Math. Soc. 330 (1992), no. 2, 545--574.

\bibitem{K83} T.~Kohno,
            {\em On the holonomy Lie algebra and the nilpotent completion of the fundamental group of the complement of hypersurfaces},
            Nagoya Math J. \textbf{92} (1983), 21--37.

\bibitem{K85} T.~Kohno,
            {\em S\'{e}rie de {P}oincar\'{e}-Koszul associ\'{e}e aux groupes de tresses pures}, 
             Invent. Math. \textbf{82} (1985), 57--75.

\bibitem{ks} J.~Kung, H.~Schenck, 
             {\em Derivation modules of orthogonal duals of hyperplane arrangements},
             Journal of Algebraic Combinatorics, {\bf 24} (2006), 253-262. 

\bibitem{LY}  A.~Libgober, S.~Yuzvinsky,
{\em Cohomology of the Orlik-Solomon algebras and local systems},
Compositio Math. \textbf{121} (2000), 337--361.

\bibitem{LS}  P.~Lima-Filho, H.~Schenck,
{\em The holonomy Lie algebra of subarrangements of $A_n$},
IMRN, (2009), 1421-1432.

\bibitem{Massey}  W.~Massey, 
{\em Completion of link modules}, Duke Math. J. 
\textbf{47} (1980), 399--420. 

\bibitem{M} J.~Morgan,
            {\em The algebraic topology on smooth algebraic varieties},
            Inst. Hautes \'Etudes Sci. Publ. Math. \textbf{48} (1978), 137-204.

\bibitem{OS} P.~Orlik, L.~Solomon,
{\em Combinatorics and topology of complements of hyperplanes}, 
Invent. Math. \textbf{56} (1980), 167--189. 

\bibitem{ot} P.~Orlik, H.~Terao,
        {\em Arrangements of Hyperplanes}, Grundlehren Math. Wiss., Bd. ~300,
        Springer-Verlag, Berlin-Heidelberg-New York, 1992.


\bibitem{othyper} P.~Orlik, H.~Terao,
        {\em Arrangements and hypergeometric integrals},
         MSJ Memoirs, \textbf{9}, Mathematical Society of Japan, Tokyo, 2007.

\bibitem{ot1} P.~Orlik, H.~Terao,
        {\em Commutative algebras for arrangements},
         Nagoya Math. J., \textbf{134} (1994), 65-73.

\bibitem{PS}  S.~Papadima, A.~Suciu,
            {\em When does the associated graded Lie algebra of an arrangement group decompose?},
            Commentarii Mathematici Helvetici, \textbf{81} (2006), 859--875.

\bibitem{PSChen} S.~Papadima, A.~Suciu,
            {\em Chen Lie algebras}, 
             IMRN \textbf{21} (2004), 1057--1086.  


\bibitem{PS2}  S.~Papadima, A.~Suciu,
            {\em Higher homotopy groups of complements of complex hyperplane arrangements}
            Adv. Math.  \textbf{165}  (2002),  71--100.

\bibitem{PY}  S.~Papadima, S.~Yuzvinsky,
            {\em On rational $K[\pi,1]$ spaces and Koszul algebras},
            J. Pure Applied Algebra \textbf{144} (1999), 157--167.


\bibitem{P} I.~Peeva,
            {\em Hyperplane arrangements and linear strands in resolutions},
             Trans. Amer. Math. Soc. \textbf{355} (2003), 609--618.

\bibitem{PerY} J.~Pereira, S.~Yuzvinsky, 
{\em Completely reducible hypersurfaces in a pencil},
Adv. Math.  \textbf{219} (2008), 672--688.
 
\bibitem{ps} N.~Proudfoot, D.~Speyer,
{\em A broken circuit ring},
Beitr\"age Algebra Geom.,  \textbf{47} (2006), 161-166.

\bibitem{R} R.~Randell, 
{\em The fundamental group of the complement of a union of complex hyperplanes},
Invent. Math. \textbf{69} (1982), 103--108.

\bibitem{Roos} J.E ~Roos,
            {\em The homotopy Lie algebra of a complex hyperplane arrangement is
 not necessarily finitely generated},
             Experimental Mathematics \textbf{17} (2008), 129--143.


\bibitem{Ry} G.~Rybnikov,
{\em On the fundamental group of the complement of a complex hyperplane arrangement},
DIMACS Tech. Report 94-13 (1994), 33--50. 

\bibitem{s} K.~Saito,
        {\em Quasihomogene isolierte Singularit\"aten von Hyperfl\"achen},
        Inventiones Mathematicae \textbf{14} (1971), 123-142.

\bibitem{slog} K. ~Saito,
        {\em Theory of logarithmic differential forms and logarithmic vector fields},
        J. Fac. Sci. Univ. Tokyo Sect. IA Math. \textbf{27} (1980), 265--291.

\bibitem{Sal}  M.~Salvetti,
{\em  Topology of the complement of real hyperplanes in $\C^N$},
Invent. Math. \textbf{88} (1987), 603--618.

\bibitem{STV} V.~Schechtman, H.~Terao, A.~Varchenko, 
{\em Local systems over complements of hyperplanes and the Kac-Kazhdan conditions for singular vectors},
J. Pure Appl. Algebra  \textbf{100} (1995), 93--102. 

\bibitem{syzRes} H.~Schenck,
{\em Resonance varieties via blowups of $\P^2$ and scrolls}, 
      IMRN, \textbf{20}, (2011), 4756-4778.

\bibitem{SS} H.~Schenck, A.~Suciu,
            {\em Lower central series and free resolutions of hyperplane arrangements},
             Trans. Amer. Math. Soc. \textbf{354} (2002), 3409--3433.

\bibitem{SS2} H.~Schenck, A.~Suciu,
{\em Resonance, linear syzygies, Chen groups, and the Bernstein-Gelfand-Gelfand correspondence},
Trans. Amer. Math. Soc. \textbf{358} (2006), 2269-2289.

\bibitem{STY} H.~Schenck, H.~Terao, M.~Yoshinaga,
{\em Logarithmic forms for quasihomogeneous curve configurations in $\P^2$},
in preparation, 2009.

\bibitem{ST1} H.~Schenck, S.~Tohaneanu,
{\em Freeness of Conic-Line arrangements in $\P^2$},
Commentarii Mathematici Helvetici, \textbf{84} (2009), 235-258.

\bibitem{ST2} H.~Schenck, S.~Tohaneanu,
{\em The Orlik-Terao algebra and 2-formality},
Math. Res. Lett.  \textbf{16}  (2009),  171-182.

\bibitem{SY} B.~Shelton, S.~Yuzvinsky,
            {\em Koszul algebras from graphs and hyperplane arrangements},
            J. London Math. Soc. \textbf{56} (1997), 477--490.

\bibitem{solt} L. ~Solomon, H. ~Terao,
        {\em A formula for the characteristic polynomial of an arrangement},
         Advances in Mathematics \textbf{64} (1987), 305-325.

\bibitem{stan} R.~Stanley,
        {\em Supersolvable lattices},
        Algebra Universalis \textbf{2} (1972), 197--217.
 
\bibitem{Su} A.~Suciu,
            {\em Fundamental groups of line arrangements: Enumerative aspects}, 
            Contemporary Math., vol.~276, Amer. Math. Soc, Providence,
            RI, 2001, pp. 43--79.

\bibitem{S} D.~Sullivan,
            {\em Infinitesimal computations in topology},
            Inst. Hautes \'Etudes Sci. Publ. Math. \textbf{47} (1977), 269--331.

\bibitem{t} H.~Terao,
        {\em Generalized exponents of a free arrangement of hyperplanes and Shepard-Todd-Brieskorn formula},
        Invent. Math. \textbf{63} (1981), 159--179.

\bibitem{t2} H.~Terao,
    {\em Arrangements of hyperplanes and their freeness I},
    J.Fac.Science Univ.Tokyo \textbf{27} (1980), 293--312.

\bibitem{T} H.~Terao,
            {\em Modular elements of lattices and topological fibration},
            Adv. Math.  \textbf{62}  (1986), 135--154.

\bibitem{tunpub} H. ~Terao,
        {\em On the homological dimensions of arrangements},
        Unpublished manuscript, 1990.

\bibitem{ter} H.~Terao,
        {\em Algebras generated by reciprocals of linear forms},
        J. Algebra,  \textbf{250} (2002), 549--558.

\bibitem{var} A.~Varchenko, 
{\em Special functions, KZ type equations, and representation theory},
        CBMS Conference Series \textbf{98}, AMS, 2003.

\bibitem{wy} M.~Wakefield, S.~Yuzvinsky,
{\em Derivations of an effective divisor on the complex projective line}, 
        Trans. Amer. Math. Soc. \textbf{359} (2007), 4389--4403.

\bibitem{y1} M.~Yoshinaga,
{\em Characterization of a free arrangement and conjecture of Edelman and Reiner},
Invent. Math. \textbf{157} (2004), no.~2, 449--454.

\bibitem{y2} M.~Yoshinaga,
{\em On the freeness of $3$-arrangements},
Bull. London Math. Soc. \textbf{37} (2005), no.~1, 126--134.

\bibitem{yu1} S.~Yuzvinsky,
        {\em First two obstructions to the freeness of arrangements},
        Transactions of the A.M.S., \textbf{335} (1993), 231-244.

\bibitem{Yuz} S.~Yuzvinsky,
{\em Cohomology of Brieskorn-Orlik-Solomon algebras},
Comm. Algebra \textbf{23} (1995), 5339--5354.

\bibitem{Yuzsurvey} S. ~Yuzvinsky,
       {\em Orlik-Solomon algebras in algebra and topology},
        Russian Math. Surveys \textbf{56} (2001), no.~2, 293--364.
 
\bibitem{yu2} S. ~Yuzvinsky,
        {\em Realization of finite Abelian groups by nets in $\P^2$},
         Compos. Math.  \textbf{140}  (2004),  1614--1624.

\bibitem{yu3} S. ~Yuzvinsky,
        {\em A new bound on the number of special fibers in a pencil of plane curves},
         Proc. Amer. Math. Soc.  \textbf{137}  (2009), 1641--1648.

\bibitem{z} G.~Ziegler,
        {\em Multiarrangements of hyperplanes and their freeness},
        Contemporary Mathematics \textbf{90}, 345-359.
        AMS, Providence, 1989.

\bibitem{z2} G. ~Ziegler,
{\em Combinatorial construction of logarithmic differential forms},
 Adv. Math. \textbf{76} (1989), 116--154.
\end{thebibliography}
\end{document}